\documentclass[12pt,a4paper]{article}

\usepackage{tensor}

\usepackage{amssymb,amsmath,amsthm,epsfig}
\numberwithin{equation}{section}
\numberwithin{figure}{section}
\numberwithin{table}{section}

\usepackage{mathrsfs}

\usepackage[table,xcdraw]{xcolor}

\usepackage{latexsym, enumerate}
\usepackage{eepic}
\usepackage{epic}
\usepackage{color}
\usepackage{ifpdf}

\usepackage[numbers,sort&compress]{natbib}
\usepackage{natbib}

\usepackage{verbatim}
\usepackage{titletoc}

\usepackage{dsfont}
\usepackage{multirow}
\usepackage{makecell}
\usepackage{bm}
\usepackage{epstopdf}
\usepackage{graphicx}
\usepackage{subfig}

\usepackage{tikz}

\graphicspath{{tuii/}{tuii/figure/}{tuii/adrIP/}{tuii/twoDimDatasets/}{tuii/locatesource/}{tuii/inverseopenarc/}{tuii/beam/}}

\usepackage{hyperref}
\hypersetup{hypertex=true,
    colorlinks=true,
    linkcolor=blue,
    filecolor=blue,
    urlcolor=blue,
    citecolor=blue,
    bookmarksopen=true,
}
\usepackage{soul}
\usepackage{xcolor}
\usepackage{color}
\sethlcolor{green}
\soulregister\cite7 
\soulregister\citep7 
\soulregister\citet7 
\soulregister\ref7 
\soulregister\pageref7 

\usepackage[toc,page]{appendix}  

\usepackage{booktabs}
\usepackage{caption}

\usepackage{algorithm}
\usepackage{algpseudocode}

\theoremstyle{plain}
\newtheorem{lem}{Lemma}[section]
\newtheorem{thm}[lem]{Theorem}
\newtheorem{cor}[lem]{Corollary}
\newtheorem{prop}[lem]{Proposition}

\theoremstyle{definition}
\newtheorem{defn}{Definition}[section]

\theoremstyle{remark}
\newtheorem{rem}{Remark}[section]

\newcommand{\diff}{\,\mathrm{d}}

\begin{document}
\title{\LARGE\bf  Geometric Optics Approximation Sampling: A Far-Field Reflector-Induced Transport Framework}

\author{
Zejun Sun\thanks
{School of Mathematics, Hunan University, Changsha 410082, China.
Email: sunzejun@hnu.edu.cn}
\and
Guang-Hui Zheng\thanks
{School of Mathematics, Hunan Provincial Key Laboratory of Intelligent Information Processing and Applied Mathematics, Hunan University, Changsha 410082, China.
Email: zhenggh2012@hnu.edu.cn (Corresponding author)}
}

\date{}
\maketitle
\begin{center}{\bf ABSTRACT}
\end{center}\smallskip
We develop a far-field geometric optics approximation sampling (GOAS) framework for constructing direct samplers from target measures. The method exploits the connection between the far-field reflector problem and optimal transport with logarithmic cost, leading to a natural primal--dual transport structure. The associated dual reflector provides a reciprocal backward transport and, in the invertible smooth setting, the inverse of the forward reflector map. For numerical realization, we adopt a supporting hyperparaboloid construction based on a discrete approximation of the target measure. This construction is gradient-free with respect to the target density and naturally accommodates both density-based and sample-based target representations. The resulting piecewise reflector admits two sampling realizations: a primal--dual consistency resampling strategy that operates directly on the nonsmooth reflector, and a softmin-regularized realization yielding an explicit smooth transport map through the physical law of reflection. We establish the well-posedness and stability of the reflector-induced sampling measure and derive Wasserstein error estimates. Numerical experiments on non-Gaussian targets and Bayesian inverse problems demonstrate the accuracy, stability, Wasserstein error behavior, and applicability of the proposed framework, and compare it with MCMC and polynomial transport-map methods.



\section{Introduction}\label{sec:intro}

Sampling from complex probability distributions is a fundamental task in computational statistics, uncertainty quantification, and Bayesian inverse problems. Markov chain Monte Carlo (MCMC) methods remain among the most widely used approaches because of their generality and asymptotic exactness \cite{robert2004,tierney1994,steve2011}. Considerable effort has been devoted to improving their efficiency through adaptive proposals \cite{haario2001,andrieu2006}, gradient- and geometry-informed dynamics such as Metropolis-adjusted Langevin algorithm and Hamiltonian Monte Carlo \cite{roberts1998,neal2011,girolami2011}, and likelihood-informed or function-space constructions for Bayesian inverse problems \cite{cotter2013,cui2016,martin2012}. Nevertheless, MCMC remains sequential and generally produces correlated samples, while sophisticated proposals may require tuning, derivative information, or additional forward and adjoint model evaluations.

Transport-based sampling provides a complementary viewpoint. 
Rather than constructing a Markov chain, one seeks a deterministic map that pushes a simple reference measure to a prescribed target measure. Once an accurate map is available, independent samples from the reference distribution can be pushed forward to generate independent samples from the target distribution.
Representative approaches include measure-transport and triangular-map methods \cite{el2012,marzouk2016,parno2018,baptista2023}, as well as normalizing flows and their continuous variants \cite{rezende2015,kobyzev2020,papamakarios2021,durkan2019,grathwohl2019,lipman2023}. These methods can be highly expressive, but typically require a prescribed map parameterization or neural-network architecture together with numerical optimization or training. Their performance therefore depends on the chosen representation and its approximation quality. Moreover, although some transport constructions are invertible by design, an inverse transport may otherwise require a separate numerical inversion or an additional approximation.

Geometric Optics Approximation Sampling (GOAS) provides a different mechanism for constructing transport. In our previous work \cite{sun2024}, sampling was connected with a near-field reflector system: a reflecting surface was constructed so that the physical law of reflection transported a simple source measure toward a prescribed target measure. The present work develops the far-field counterpart, but the distinction is not merely geometric. The near-field reflector problem is not, in general, a classical optimal transport problem; its transport cost depends nonlinearly on the reflector potential, and weak solutions need not minimize the corresponding transport functional \cite{karakhanyan2010,graf2012}. In contrast, the far-field reflector antenna problem admits a classical optimal transport formulation with logarithmic cost \cite{wang2004}. This additional structure is the central feature exploited in this paper.

More specifically, a far-field reflector induces a transport from an incident direction on the source sphere to a reflected direction on the output sphere. After a diffeomorphic transformation between the output sphere and the Euclidean target domain, this reflector-induced transport becomes a direct sampling mechanism. The optimal transport structure further provides a natural primal--dual reflector geometry. If $R_\rho$ denotes the primal reflector and $R_{\rho^\ast}$ its dual reflector, their associated subdifferential maps satisfy the reciprocal relation
\[
y\in T_\rho(x)
\quad\Longleftrightarrow\quad
x\in T_{\rho^\ast}(y)
\]
almost everywhere. Thus, the dual reflector supplies a geometric inverse relation to the forward reflector transport. 
This primal--dual structure is specific to the far-field formulation, leading naturally to reciprocal forward and inverse reflector-induced transports.

The numerical framework developed in this paper follows this structure. We adopt the supporting hyperparaboloid method as a constructive realization of the far-field reflector problem.
We construct a discrete far-field reflector by extending the classical supporting paraboloid method \cite{caffarelli1999,kochengin2003} to a supporting hyperparaboloid formulation in arbitrary dimension. The construction is gradient-free with respect to the target density and only requires a discrete approximation of the target measure. Since the resulting reflector is generally continuous but nonsmooth at interfaces between supporting paraboloid patches, we introduce a primal--dual consistency resampling strategy that exploits the reciprocal relation between the primal and dual reflectors to generate geometrically consistent samples.

Alternatively, we smooth the piecewise reflector by a softmin regularization. 
The resulting smooth surface induces a transport map explicitly through the physical law of reflection, so that sampling requires only the evaluation of this reflector map. 
In the far-field setting, the corresponding dual reflector further provides the inverse transport map.
Although only the forward map is required for direct sampling in the present work, the inverse transport map arises naturally as a byproduct of the primal--dual geometry and may be useful in other transport-based constructions, such as transport-informed MCMC proposals.

The theoretical analysis focuses on the measure induced by the far-field reflector transport. 
We first define the geometric optics approximation measure associated with the reflector; in the single-valued setting, it is precisely the push-forward of the source measure under the composition of the reflector-induced transport map and a projection.
Its well-posedness and stability provide the mathematical basis for interpreting the reflector map as a transport mechanism for direct sampling.
We then quantify the approximation error of the measure in Wasserstein distance, accounting for target discretization, numerical evaluation of the induced masses, reflector construction tolerance, and softmin smoothing.

The main contributions of this work are summarized as follows.
\begin{itemize}

\item
We develop a far-field formulation of GOAS that places reflector-induced sampling within a classical optimal transport framework. In contrast to the near-field formulation, the far-field reflector problem admits an optimal transport structure with logarithmic cost. The resulting transport is determined by the physical law of reflection, rather than by a prescribed parametric map or a trained neural network.

\item
We exploit the intrinsic primal--dual structure of the far-field reflector problem. The dual reflector is obtained naturally from reflector duality and induces the reciprocal transport relation. Consequently, the far-field formulation provides not only the forward transport from the source to the target, but also a natural mechanism for inverse transport. This primal--dual structure also underlies the consistency relation used in the sampling realization for nonsmooth reflectors.

\item
We develop a gradient-free numerical realization based on supporting hyperparaboloid in general dimensions. The method constructs a piecewise reflector from a discrete approximation of the target measure, thereby naturally accommodating both density-based and sample-based target representations. Since the construction uses only target masses or density evaluations at prescribed target points, no gradients of the target density are required. The resulting piecewise reflector admits two complementary sampling realizations: primal--dual consistency resampling can be applied directly to the nonsmooth reflector, while softmin regularization produces a smooth reflector from which an explicit transport map is obtained through the reflection law. The supporting hyperparaboloid method is one numerical realization of the far-field GOAS framework, and other solvers for the corresponding reflector or optimal transport problem may also be used.

\item
We establish a measure-theoretic and quantitative analysis of the resulting sampling framework. 
We define the geometric optics approximation measure and study its well-posedness and stability, providing the mathematical basis for the far-field GOAS.
We further derive Wasserstein error estimates that account for target discretization, numerical evaluation of the induced masses, reflector construction tolerance, and softmin smoothing. 
Numerical experiments on non-Gaussian distributions and Bayesian inverse problems illustrate the performance of the proposed framework and compare it with MCMC and transport-map-based methods.
\end{itemize}

The remainder of the paper is organized as follows. Section~\ref{sec:GOASmethod} develops the far-field reflector geometry, its optimal transport formulation, the primal--dual reflector structure, the geometric optics approximation measure, and the associated well-posedness results. Section~\ref{sec:GOASAlg} presents the numerical realization of far-field GOAS, including target discretization, supporting hyperparaboloid construction, primal--dual consistency resampling, softmin smoothing, and the complete sampling algorithm. Section~\ref{sec:wasserstein_error} derives Wasserstein error estimates for the numerical GOAS measure. Section~\ref{sec:numexp} presents some numerical experiments. Section~\ref{sec:conclusion} concludes the paper and discusses future directions.
Additional geometric formulation of far-field reflector, proofs, and numerical details are provided in the appendices.

\section{Far-field reflector-induced optimal transport}\label{sec:GOASmethod}

This section develops the mathematical structure underlying far-field GOAS. We first recall the geometry of convex far-field reflectors and their subdifferential maps. We then formulate the reflector problem as an optimal transport problem with logarithmic cost and introduce the dual reflector, which provides the reciprocal transport relation used later in the primal--dual sampler. Finally, we define the geometric optics approximation measure and recall the existence and stability results that provide the measure-theoretic foundation of the method.

\begin{figure}[htbp]
  \centering


\tikzset {_hbwj7gujh/.code = {\pgfsetadditionalshadetransform{ \pgftransformshift{\pgfpoint{89.1 bp } { -128.7 bp }  }  \pgftransformscale{1.32 }  }}}
\pgfdeclareradialshading{_1io0q7yh0}{\pgfpoint{-72bp}{104bp}}{rgb(0bp)=(1,1,1);
rgb(0bp)=(1,1,1);
rgb(25bp)=(0.99,0.02,0.02);
rgb(400bp)=(0.99,0.02,0.02)}


\tikzset{
pattern size/.store in=\mcSize,
pattern size = 5pt,
pattern thickness/.store in=\mcThickness,
pattern thickness = 0.3pt,
pattern radius/.store in=\mcRadius,
pattern radius = 1pt}
\makeatletter
\pgfutil@ifundefined{pgf@pattern@name@_rvuk3ivf4}{
\pgfdeclarepatternformonly[\mcThickness,\mcSize]{_rvuk3ivf4}
{\pgfqpoint{0pt}{-\mcThickness}}
{\pgfpoint{\mcSize}{\mcSize}}
{\pgfpoint{\mcSize}{\mcSize}}
{
\pgfsetcolor{\tikz@pattern@color}
\pgfsetlinewidth{\mcThickness}
\pgfpathmoveto{\pgfqpoint{0pt}{\mcSize}}
\pgfpathlineto{\pgfpoint{\mcSize+\mcThickness}{-\mcThickness}}
\pgfusepath{stroke}
}}
\makeatother


\tikzset{
pattern size/.store in=\mcSize,
pattern size = 5pt,
pattern thickness/.store in=\mcThickness,
pattern thickness = 0.3pt,
pattern radius/.store in=\mcRadius,
pattern radius = 1pt}
\makeatletter
\pgfutil@ifundefined{pgf@pattern@name@_4obmmm7ro}{
\pgfdeclarepatternformonly[\mcThickness,\mcSize]{_4obmmm7ro}
{\pgfqpoint{0pt}{-\mcThickness}}
{\pgfpoint{\mcSize}{\mcSize}}
{\pgfpoint{\mcSize}{\mcSize}}
{
\pgfsetcolor{\tikz@pattern@color}
\pgfsetlinewidth{\mcThickness}
\pgfpathmoveto{\pgfqpoint{0pt}{\mcSize}}
\pgfpathlineto{\pgfpoint{\mcSize+\mcThickness}{-\mcThickness}}
\pgfusepath{stroke}
}}
\makeatother

\tikzset{every picture/.style={line width=0.75pt}} 

\begin{tikzpicture}[x=0.75pt,y=0.75pt,yscale=-1,xscale=1]

\draw [fill={rgb, 255:red, 242; green, 236; blue, 217 }  ,fill opacity=1 ]   (231.14,94.04) .. controls (192.01,98.53) and (123.99,186.19) .. (128.48,239.32) .. controls (132.97,292.46) and (202.91,321.86) .. (247.98,320.42) ;
\draw  [fill={rgb, 255:red, 242; green, 236; blue, 217 }  ,fill opacity=1 ] (210.06,116.06) .. controls (193.55,153.48) and (194.09,224.4) .. (211.25,274.48) .. controls (228.42,324.55) and (255.71,334.81) .. (272.21,297.4) .. controls (288.71,259.98) and (288.18,189.05) .. (271.02,138.98) .. controls (253.85,88.9) and (226.56,78.64) .. (210.06,116.06) -- cycle ;
\draw [draw opacity=0][fill={rgb, 255:red, 253; green, 198; blue, 204 }  ,fill opacity=1 ]   (260.86,141.83) .. controls (242.25,157.82) and (236.62,198.96) .. (237.85,211.8) .. controls (239.08,224.64) and (249.57,265.48) .. (270.57,274.78) ;
\draw  [fill={rgb, 255:red, 253; green, 198; blue, 204 }  ,fill opacity=1 ] (256.21,208.58) .. controls (254.03,171.69) and (256.92,141.53) .. (262.66,141.21) .. controls (268.41,140.89) and (274.83,170.52) .. (277.02,207.41) .. controls (279.2,244.29) and (276.32,274.45) .. (270.57,274.78) .. controls (264.83,275.1) and (258.4,245.46) .. (256.21,208.58) -- cycle ;

\draw (412.66,207.25) node  {\includegraphics[width=56.5pt,height=89.89pt]{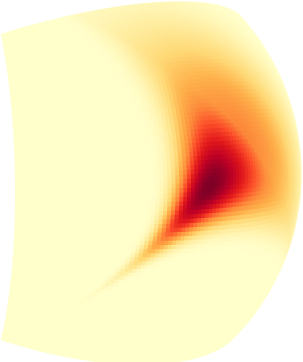}};
\draw   (239.07,203.61) .. controls (239.07,146.94) and (286.29,101) .. (344.53,101) .. controls (402.78,101) and (450,146.94) .. (450,203.61) .. controls (450,260.28) and (402.78,306.22) .. (344.53,306.22) .. controls (286.29,306.22) and (239.07,260.28) .. (239.07,203.61) -- cycle ;
\draw  [draw opacity=0][fill={rgb, 255:red, 254; green, 83; blue, 104 }  ,fill opacity=1 ] (338.72,203.61) .. controls (338.72,200.87) and (341.32,198.64) .. (344.53,198.64) .. controls (347.75,198.64) and (350.35,200.87) .. (350.35,203.61) .. controls (350.35,206.36) and (347.75,208.58) .. (344.53,208.58) .. controls (341.32,208.58) and (338.72,206.36) .. (338.72,203.61) -- cycle ;
\draw (549.27,292.34) node [xslant=0.54,xscale=-1] {\includegraphics[width=102.05pt,height=43.61pt]{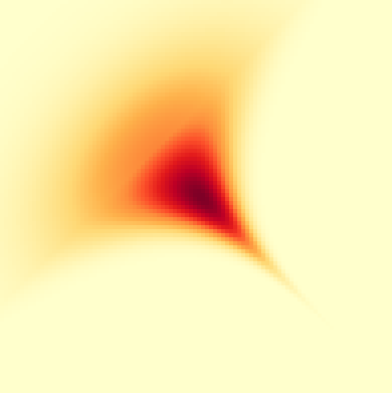}};
\draw  [draw opacity=0] (239.7,198.96) .. controls (239.48,200.5) and (239.36,202.05) .. (239.36,203.61) .. controls (239.36,235.54) and (286.45,261.43) .. (344.53,261.43) .. controls (401.84,261.43) and (448.44,236.23) .. (449.68,204.9) -- (344.53,203.61) -- cycle ; \draw   (239.7,198.96) .. controls (239.48,200.5) and (239.36,202.05) .. (239.36,203.61) .. controls (239.36,235.54) and (286.45,261.43) .. (344.53,261.43) .. controls (401.84,261.43) and (448.44,236.23) .. (449.68,204.9) ;
\draw  [draw opacity=0][dash pattern={on 4.5pt off 4.5pt}] (240.06,203.28) .. controls (239.6,201.76) and (239.36,200.21) .. (239.36,198.64) .. controls (239.36,176.47) and (286.45,158.5) .. (344.53,158.5) .. controls (401.49,158.5) and (447.88,175.78) .. (449.65,197.36) -- (344.53,198.64) -- cycle ; \draw  [dash pattern={on 4.5pt off 4.5pt}] (240.06,203.28) .. controls (239.6,201.76) and (239.36,200.21) .. (239.36,198.64) .. controls (239.36,176.47) and (286.45,158.5) .. (344.53,158.5) .. controls (401.49,158.5) and (447.88,175.78) .. (449.65,197.36) ;
\draw [color={rgb, 255:red, 250; green, 141; blue, 141 }  ,draw opacity=1 ]   (344.53,203.61) -- (133.15,257.99) -- (210.99,213.87) ;
\draw [shift={(212.73,212.88)}, rotate = 150.46] [color={rgb, 255:red, 250; green, 141; blue, 141 }  ,draw opacity=1 ][line width=0.75]    (10.93,-3.29) .. controls (6.95,-1.4) and (3.31,-0.3) .. (0,0) .. controls (3.31,0.3) and (6.95,1.4) .. (10.93,3.29)   ;
\draw    (559.46,273.23) .. controls (519.6,290.26) and (485.14,159.65) .. (419.02,165.39) ;
\draw [shift={(418.02,165.49)}, rotate = 353.88] [color={rgb, 255:red, 0; green, 0; blue, 0 }  ][line width=0.75]    (10.93,-3.29) .. controls (6.95,-1.4) and (3.31,-0.3) .. (0,0) .. controls (3.31,0.3) and (6.95,1.4) .. (10.93,3.29)   ;
\draw    (396.78,223.55) .. controls (450.26,212.95) and (474,328.83) .. (531.01,308.23) ;
\draw [shift={(532.75,307.56)}, rotate = 157.93] [color={rgb, 255:red, 0; green, 0; blue, 0 }  ][line width=0.75]    (10.93,-3.29) .. controls (6.95,-1.4) and (3.31,-0.3) .. (0,0) .. controls (3.31,0.3) and (6.95,1.4) .. (10.93,3.29)   ;
\draw  [color={rgb, 255:red, 254; green, 127; blue, 127 }  ,draw opacity=1 ] (288.86,213.56) -- (278.29,221.07) -- (291.19,222.38) ;
\draw [color={rgb, 255:red, 218; green, 183; blue, 124 }  ,draw opacity=1 ]   (133.15,257.99) -- (76.99,280.15) ;
\draw [shift={(75.13,280.88)}, rotate = 338.47] [color={rgb, 255:red, 218; green, 183; blue, 124 }  ,draw opacity=1 ][line width=0.75]    (10.93,-3.29) .. controls (6.95,-1.4) and (3.31,-0.3) .. (0,0) .. controls (3.31,0.3) and (6.95,1.4) .. (10.93,3.29)   ;

\draw (341.04,275.32) node [anchor=north west][inner sep=0.75pt]   [align=left] {$\displaystyle S^{2}$};
\draw (445.53,272.88) node [anchor=north west][inner sep=0.75pt]   [align=left] {$\displaystyle Q$};
\draw (496.26,194.9) node [anchor=north west][inner sep=0.75pt]   [align=left] {$\displaystyle Q^{-1}$};
\draw (252.45,233.43) node [anchor=north west][inner sep=0.75pt]   [align=left] {$\displaystyle \Omega _{I}$};
\draw (384.51,227.36) node [anchor=north west][inner sep=0.75pt]   [align=left] {$\displaystyle \Omega _{O}$};
\draw (580.76,300.18) node [anchor=north west][inner sep=0.75pt]   [align=left] {$\displaystyle \Omega $};
\draw (592.82,265.49) node [anchor=north west][inner sep=0.75pt]   [align=left] {$\displaystyle \mu _{t}$};
\draw (386.51,166.63) node [anchor=north west][inner sep=0.75pt]   [align=left] {$\displaystyle \mu _{O}$};
\draw (257.46,186.21) node [anchor=north west][inner sep=0.75pt]   [align=left] {$\displaystyle \mu _{I}$};
\draw (338.91,179.65) node [anchor=north west][inner sep=0.75pt]   [align=left] {$\displaystyle \mathcal{O}$};
\draw (288.23,196.06) node [anchor=north west][inner sep=0.75pt]   [align=left] {$\displaystyle x$};
\draw (182.39,202.83) node [anchor=north west][inner sep=0.75pt]   [align=left] {$\displaystyle y$};
\draw (170.18,249.13) node [anchor=north west][inner sep=0.75pt]   [align=left] {$\displaystyle \rho ( x)$};
\draw (216.78,134.95) node [anchor=north west][inner sep=0.75pt]   [align=left] {$\displaystyle R_{\rho }$};
\draw (84.73,251.32) node [anchor=north west][inner sep=0.75pt]   [align=left] {$\displaystyle \nu$};

\end{tikzpicture}
\caption{Schematic of far-field GOAS. A source direction $x\sim\mu_I$ is transported by the reflector-induced map to an output direction $y=T_\rho(x)$ and then mapped to the target domain by $z=Q(y)$.}
\label{fig:RefAntennaSystem}
\end{figure}

\subsection{Reflector-induced optimal transport}\label{sec:farfield_OT}
Let $\mathcal O$ be a point source located at the origin of $\mathbb R^{n+1}$. The incident directions are contained in an input domain $\Omega_I\subset S^n_+$, where
$
S^n_+=\{x\in S^n:x_{n+1}>0\},
$
and reflected rays propagate toward an output domain $\Omega_O\subset S^n_-$, with
$
S^n_-=\{y\in S^n:y_{n+1}<0\},
$
where $S^n$ is the unit sphere.
A reflector is represented as the radial graph
\begin{align}\label{RGraph}
R_\rho=\{x\rho(x):x\in\Omega_I\},
\end{align}
where $\rho:\Omega_I\to(0,\infty)$ is the polar radius. The law of reflection sends an incident direction $x\in\Omega_I$ to
\begin{align}\label{RefDir}
y=T(x)=T_\rho(x)=x-2(x\cdot\nu)\nu,
\end{align}
where $\nu$ is the outward unit normal of the reflector. Thus, the reflector geometry determines the reflection map $T$ explicitly.
Let the target measure on a domain $\Omega\subset\mathbb R^n$ be
\begin{align}\label{TarMea}
\mu_t(U)=\int_U\pi(z)\diff\mu(z),
\qquad U\subset\Omega,
\end{align}
where $\mu$ denotes Lebesgue measure. Let $Q$ be a diffeomorphic transformation from the relevant part of the sphere to the target domain. The target measure is represented on the output sphere by
\begin{align}\label{TarMeaSph}
\mu_O(w)
=
\int_w \pi\circ Q(y)|J(Q)|\diff\sigma(y)
=
\int_w g(y)\diff\sigma(y),
\qquad w\subset\Omega_O,
\end{align}
where $\sigma$ denotes spherical measure. The source measure is
\begin{align}\label{eq:source_measure}
\mu_I(E)=\int_E f(x)\diff\sigma(x),
\qquad E\subset\Omega_I.
\end{align}

Then once the reflecting surface has been constructed, samples from the target measure are generated by first drawing $x\sim\mu_I$, reflecting it to $y=T_\rho(x)$, and then mapping the reflected direction to the target domain:
\begin{align}\label{MapSam}
z=Q\circ T_\rho(x).
\end{align}
This relation is the sampling interpretation of far-field reflector-induced optimal transport.
See Figure \ref{fig:RefAntennaSystem} for sampling the non-Gaussian target distribution using the reflector-induced transport map.

Therefore, the sampling problem is reduced to the construction of an appropriate reflecting surface.
We assume the energy conservation condition
\begin{align}\label{EneCon}
\int_{\Omega_I}f\diff\sigma
=
\int_{\Omega_O}g\diff\sigma.
\end{align}
The reflector surface that transports $\mu_I$ to $\mu_O$ satisfies the Monge--Amp\`ere type equation \cite{oliker1993,wang1996}
\begin{align}\label{RhoPDE}
\frac{
\det(
-\nabla_i\nabla_j\rho
+
2\rho^{-1}\nabla_i\rho\nabla_j\rho
+
(|\nabla\rho|^2-\rho^2)\delta_{ij}/2\rho
)
}{
(|\nabla\rho|^2+\rho^2)/2\rho
}
=
\frac{f(x)}{g\circ T_\rho(x)},
\qquad x\in\Omega_I,
\end{align}
subject to the transport boundary condition
\begin{align}\label{RhoPDECon}
T_\rho(\Omega_I)=\Omega_O.
\end{align}
Here $\nabla$ denotes the covariant derivative on the sphere and $\delta_{ij}$ is the Kronecker symbol.

The essential distinction of the far-field reflector problem is that this geometric transport is also a classical optimal transport problem.

\begin{thm}\label{OTRef}
A solution of the reflector antenna problem \eqref{RhoPDE}--\eqref{RhoPDECon} is equivalent to an optimal transport map solving
\[
\inf_{T_\sharp\mu_I=\mu_O}
\int_{\Omega_I} c(x,T(x))\diff\mu_I(x),
\]
with logarithmic cost
$
c(x,y)=-\log(1-x\cdot y).
$
Here $T_\sharp\mu_I=\mu_O$ means that
$
\mu_O(w)=\mu_I(T^{-1}(w)), w\subset\Omega_O.
$
\end{thm}

The proof can be found in \cite{wang2004}. The optimal transport map is exactly the reflection map \eqref{RefDir}, which is determined explicitly by physical law of reflection.

Theorem~\ref{OTRef} places far-field GOAS directly within optimal transport: the reflector-induced map solves a Monge problem with the logarithmic cost $c(x,y)=-\log(1-x\cdot y)$. This structure is not available in the same classical form for the near-field reflector problem and is the source of the primal--dual construction developed next.

\subsubsection{Dual reflectors and inverse transport}\label{sec:dual_transport}

The optimal transport structure of the far-field reflector is accompanied by a natural dual reflector. This duality is not introduced as a separate numerical device; it is intrinsic to the reflector geometry and provides the reciprocal relation between incident and reflected directions.

\begin{figure}[htbp]
  \centering

\tikzset {_hbwj7gujh/.code = {\pgfsetadditionalshadetransform{ \pgftransformshift{\pgfpoint{89.1 bp } { -128.7 bp }  }  \pgftransformscale{1.32 }  }}}
\pgfdeclareradialshading{_1io0q7yh0}{\pgfpoint{-72bp}{104bp}}{rgb(0bp)=(1,1,1);
rgb(0bp)=(1,1,1);
rgb(25bp)=(0.99,0.02,0.02);
rgb(400bp)=(0.99,0.02,0.02)}


\tikzset{
pattern size/.store in=\mcSize,
pattern size = 5pt,
pattern thickness/.store in=\mcThickness,
pattern thickness = 0.3pt,
pattern radius/.store in=\mcRadius,
pattern radius = 1pt}
\makeatletter
\pgfutil@ifundefined{pgf@pattern@name@_rvuk3ivf4}{
\pgfdeclarepatternformonly[\mcThickness,\mcSize]{_rvuk3ivf4}
{\pgfqpoint{0pt}{-\mcThickness}}
{\pgfpoint{\mcSize}{\mcSize}}
{\pgfpoint{\mcSize}{\mcSize}}
{
\pgfsetcolor{\tikz@pattern@color}
\pgfsetlinewidth{\mcThickness}
\pgfpathmoveto{\pgfqpoint{0pt}{\mcSize}}
\pgfpathlineto{\pgfpoint{\mcSize+\mcThickness}{-\mcThickness}}
\pgfusepath{stroke}
}}
\makeatother


\tikzset{
pattern size/.store in=\mcSize,
pattern size = 5pt,
pattern thickness/.store in=\mcThickness,
pattern thickness = 0.3pt,
pattern radius/.store in=\mcRadius,
pattern radius = 1pt}
\makeatletter
\pgfutil@ifundefined{pgf@pattern@name@_4obmmm7ro}{
\pgfdeclarepatternformonly[\mcThickness,\mcSize]{_4obmmm7ro}
{\pgfqpoint{0pt}{-\mcThickness}}
{\pgfpoint{\mcSize}{\mcSize}}
{\pgfpoint{\mcSize}{\mcSize}}
{
\pgfsetcolor{\tikz@pattern@color}
\pgfsetlinewidth{\mcThickness}
\pgfpathmoveto{\pgfqpoint{0pt}{\mcSize}}
\pgfpathlineto{\pgfpoint{\mcSize+\mcThickness}{-\mcThickness}}
\pgfusepath{stroke}
}}
\makeatother
\tikzset{every picture/.style={line width=0.75pt}} 
\begin{tikzpicture}[x=0.75pt,y=0.75pt,yscale=-1,xscale=1]

\draw  [dash pattern={on 0.84pt off 2.51pt}] (295.95,224.39) .. controls (295.95,193.8) and (320.74,169) .. (351.33,169) .. controls (381.92,169) and (406.72,193.8) .. (406.72,224.39) .. controls (406.72,254.98) and (381.92,279.78) .. (351.33,279.78) .. controls (320.74,279.78) and (295.95,254.98) .. (295.95,224.39) -- cycle ;
\draw  [draw opacity=0][line width=1.5]  (280.86,308.92) .. controls (231.01,299.65) and (194,265.28) .. (194,224.29) .. controls (194,183.42) and (230.8,149.13) .. (280.43,139.74) -- (307.84,224.29) -- cycle ; \draw  [color={rgb, 255:red, 189; green, 16; blue, 224 }  ,draw opacity=1 ][line width=1.5]  (280.86,308.92) .. controls (231.01,299.65) and (194,265.28) .. (194,224.29) .. controls (194,183.42) and (230.8,149.13) .. (280.43,139.74) ;
\draw  [draw opacity=0][fill={rgb, 255:red, 255; green, 170; blue, 223 }  ,fill opacity=1 ] (343.26,224.39) .. controls (343.26,219.93) and (346.87,216.31) .. (351.33,216.31) .. controls (355.8,216.31) and (359.41,219.93) .. (359.41,224.39) .. controls (359.41,228.85) and (355.8,232.47) .. (351.33,232.47) .. controls (346.87,232.47) and (343.26,228.85) .. (343.26,224.39) -- cycle ;
\draw [color={rgb, 255:red, 250; green, 88; blue, 107 }  ,draw opacity=1 ][line width=0.75]    (350,224.53) -- (488.4,141.72) -- (340.19,165.9) ;
\draw [shift={(338.22,166.23)}, rotate = 350.73] [color={rgb, 255:red, 250; green, 88; blue, 107 }  ,draw opacity=1 ][line width=0.75]    (10.93,-3.29) .. controls (6.95,-1.4) and (3.31,-0.3) .. (0,0) .. controls (3.31,0.3) and (6.95,1.4) .. (10.93,3.29)   ;
\draw [color={rgb, 255:red, 189; green, 16; blue, 224 }  ,draw opacity=1 ][line width=0.75]    (350,224.53) -- (199.81,249.04) -- (336.5,167.25) ;
\draw [shift={(338.22,166.23)}, rotate = 149.11] [color={rgb, 255:red, 189; green, 16; blue, 224 }  ,draw opacity=1 ][line width=0.75]    (10.93,-3.29) .. controls (6.95,-1.4) and (3.31,-0.3) .. (0,0) .. controls (3.31,0.3) and (6.95,1.4) .. (10.93,3.29)   ;
\draw  [color={rgb, 255:red, 189; green, 16; blue, 224 }  ,draw opacity=1 ] (305.6,228.67) -- (297.31,233.16) -- (306.6,233.75) ;
\draw  [color={rgb, 255:red, 189; green, 16; blue, 224 }  ,draw opacity=1 ] (230.84,227.02) -- (240.3,224.96) -- (235.39,230.95) ;
\draw  [color={rgb, 255:red, 247; green, 84; blue, 104 }  ,draw opacity=1 ] (389.33,197.54) -- (398.62,195.07) -- (394.27,201.24) ;
\draw  [color={rgb, 255:red, 247; green, 84; blue, 104 }  ,draw opacity=1 ] (439.55,146.75) -- (431.8,151.05) -- (441.68,151.58) ;
\draw  [draw opacity=0][line width=1.5]  (437.46,340.54) .. controls (484.13,330.11) and (519.39,282.32) .. (519.39,224.95) .. controls (519.39,167.79) and (484.39,120.15) .. (437.98,109.48) -- (418.28,224.95) -- cycle ; \draw  [color={rgb, 255:red, 254; green, 107; blue, 125 }  ,draw opacity=1 ][line width=1.5]  (437.46,340.54) .. controls (484.13,330.11) and (519.39,282.32) .. (519.39,224.95) .. controls (519.39,167.79) and (484.39,120.15) .. (437.98,109.48) ;
\draw  [draw opacity=0] (316.26,267.77) .. controls (304.12,257.62) and (296.33,241.96) .. (296.33,224.39) .. controls (296.33,207.23) and (303.75,191.9) .. (315.41,181.74) -- (349.03,224.39) -- cycle ; \draw  [color={rgb, 255:red, 189; green, 16; blue, 224 }  ,draw opacity=1 ] (316.26,267.77) .. controls (304.12,257.62) and (296.33,241.96) .. (296.33,224.39) .. controls (296.33,207.23) and (303.75,191.9) .. (315.41,181.74) ;
\draw  [draw opacity=0] (384.75,269.46) .. controls (398.06,259.46) and (406.72,243.15) .. (406.72,224.72) .. controls (406.72,206.94) and (398.65,191.12) .. (386.13,181.06) -- (354.03,224.72) -- cycle ; \draw  [color={rgb, 255:red, 250; green, 98; blue, 98 }  ,draw opacity=1 ] (384.75,269.46) .. controls (398.06,259.46) and (406.72,243.15) .. (406.72,224.72) .. controls (406.72,206.94) and (398.65,191.12) .. (386.13,181.06) ;

\draw (360.03,223.22) node [anchor=north west][inner sep=0.75pt]   [align=left] {$\displaystyle \mathcal{O}$};
\draw (411.22,210.39) node [anchor=north west][inner sep=0.75pt]   [align=left] {$\displaystyle \Omega _{O}$};
\draw (273.33,210.33) node [anchor=north west][inner sep=0.75pt]   [align=left] {$\displaystyle \Omega _{I}$};
\draw (523.33,209.17) node [anchor=north west][inner sep=0.75pt]   [align=left] {$\displaystyle R_{\rho ^{\ast }}$};
\draw (169.67,210.17) node [anchor=north west][inner sep=0.75pt]   [align=left] {$\displaystyle R_{\rho }$};
\draw (344.45,255.77) node [anchor=north west][inner sep=0.75pt]   [align=left] {$\displaystyle S$};
\draw (226.83,206.5) node [anchor=north west][inner sep=0.75pt]   [align=left] {$\displaystyle y$};
\draw (375.17,184.04) node [anchor=north west][inner sep=0.75pt]   [align=left] {$\displaystyle y$};
\draw (418.33,134.85) node [anchor=north west][inner sep=0.75pt]   [align=left] {$\displaystyle x$};
\draw (303.5,213.8) node [anchor=north west][inner sep=0.75pt]   [align=left] {$\displaystyle x$};
\draw (449,164.68) node [anchor=north west][inner sep=0.75pt]   [align=left] {$\displaystyle \rho ^{\ast }( y)$};
\draw (230.67,243.78) node [anchor=north west][inner sep=0.75pt]   [align=left] {$\displaystyle \rho ( x)$};
\end{tikzpicture}
\caption{The primal–dual geometry structure}
\label{fig:dualRef}
\end{figure}

The dual reflector is defined through the Legendre-type transform.
Let $R_{\rho}$ be a convex reflector surface.
The Legendre transform of $\rho$ with respect to the $1/(1-x\cdot y),x\in\Omega_I,y\in\Omega_O$ is a function $\rho^\ast$ defined on $S^n$, given by
\begin{align}\label{rhoLT}
\rho^\ast(y)=\inf_{x\in\Omega_I}\frac{1}{\rho(x)(1-x\cdot y)}.
\end{align}
The following duality relation holds.
\begin{prop}\label{LTDual}
Let $\rho^\ast$ be defined by \eqref{rhoLT}.
Then
\[
y\in T_\rho(x)
\quad\Longleftrightarrow\quad
x\in T_{\rho^\ast}(y)
\]
almost everywhere, where $x\in\Omega_I,y\in\Omega_O$, $T_{\rho}$ and $T_{\rho^\ast}$ are the subdifferential reflector maps given by Definition \ref{DefSubdiff} with respect to $\rho$ and $\rho^\ast$, respectively.
\end{prop}
The proof of this proposition has been given in \cite{wang2004}.
Thus, the dual reflector naturally provides an inverse transport mechanism associated with the primal reflector.

When the primal and dual reflectors are smooth and their associated reflector maps are single-valued and one-to-one, the set-valued relation in Proposition~\ref{LTDual} reduces to the usual inverse relation
\[
T_{\rho^\ast}(T_\rho(x))=x,
\qquad
T_\rho(T_{\rho^\ast}(y))=y,
\]
see Figure \ref{fig:dualRef}.
Hence the primal and dual reflectors naturally induce a forward--inverse pair of geometric transports. In the discrete reflector setting, the same duality will be used to construct a primal--dual consistency resampling in Section \ref{sec:primal_dual_sampling}.
In the smooth reflector setting, the same duality directly yields the inverse of the forward transport map.

\subsection{Geometric optics approximation measure}\label{sec:GOA_measure}
Using the visibility map introduced in Section~\ref{sec:farfield_geometry}, we define the measure induced by a reflector $R$ on $\Omega_O$ by
\begin{align}\label{pullbackMea}
\mu_R(B)
:=
\int_{V(B)} f(x)\,\diff\sigma(x),
\qquad
B\in\mathcal B(\Omega_O),
\end{align}
where $f\ge0$, $f\in L^1(\Omega_I)$ is the source illumination density and $\mathcal B(\Omega_O)$ denotes the Borel $\sigma$-algebra on $\Omega_O$. Under the standard properties of convex far-field reflectors, the visibility construction defines a finite Borel measure on $\Omega_O$. The quantity $\mu_R(B)$ represents the total source energy reflected into the set of directions $B$.

When the reflector is smooth, the reflector relation reduces to a single-valued measurable reflection map $T:\Omega_I\to\Omega_O$, and \eqref{pullbackMea} coincides with the usual push-forward representation,
\[
\mu_R=T_\sharp\mu_I,
\qquad
\mu_R(B)=\mu_I\bigl(T^{-1}(B)\bigr).
\]
The reflector-induced sampling measure 
\[
\mu_{\mathrm{GOAS}}
=
Q_\sharp\mu_R
\]
is referred to as the geometric optics approximation sampling (GOAS) measure.

\begin{defn}[Weak solution]\label{GOAMea}
A reflecting surface $R$ is called a weak solution of the reflector antenna problem if
\begin{align}\label{WeaSolCon}
\mu_O(w)=\mu_R(w),
\qquad \forall w\subset\Omega_O.
\end{align}
\end{defn}

Condition \eqref{WeaSolCon} is precisely the conservation of energy.
When $R$ is smooth, this weak formulation is equivalent to the reflector equation \eqref{RhoPDE} together with the boundary condition \eqref{RhoPDECon}.

\subsubsection{Well-posedness}\label{sec:wellposedness}

The GOAS measure depends on the weak solution of the reflector antenna problem.
Therefore, its well-posedness follows from the existence and uniqueness theory of reflector antenna solutions.
These results have been studied from both PDE and optimal transport perspectives
\cite{caffarelli1999,caffarelli2008,wang1996,wang2004,guan1998}.
We recall the following standard result.

\begin{thm}[Existence and uniqueness]\label{ExiUniThe}
Let $\mu_I$ and $\mu_O$ be measures supported on compact domains $\Omega_I,\Omega_O\subset S^n$, respectively.
Assume that the energy conservation condition \eqref{EneCon} holds.
Then the far-field reflector problem admits a weak solution $R_\rho$ satisfying \eqref{WeaSolCon}.
Moreover, the radial function $\rho$ is unique up to multiplication by a positive constant.
\end{thm}

This result follows from the classical existence and uniqueness theory for the far-field reflector problem; see \cite{wang1996,caffarelli1999,caffarelli2008,wang2004}.

\begin{thm}[Stability with respect to the target measure]
\label{StaTarDomThe}
Let $\{\mu_t^k\}_{k\ge1}$ and $\mu_t$ be target measures such that
\[
\mu_t^k \rightharpoonup \mu_t
\qquad\text{as }k\to\infty.
\]
Assume that the corresponding output measures
\[
\mu_O^k=(Q^{-1})_\sharp\mu_t^k,
\qquad
\mu_O=(Q^{-1})_\sharp\mu_t
\]
are supported in a common compact subset of $\Omega_O$.
Let $R_{\rho_k}$ and $R_\rho$ be the corresponding weak reflectors, normalized by
$
\rho_k(x_0)=\rho(x_0)=1
$
for some fixed $x_0\in\Omega_I$.
Then
\[
\rho_k\to\rho
\qquad\text{uniformly on }\Omega_I.
\]
Consequently,
\[
d_H(R_{\rho_k},R_\rho)\longrightarrow0
\qquad\text{as }k\to\infty,
\]
where $d_{H}(\cdot,\cdot)$ denotes the Hausdorff metric.
\end{thm}

\begin{proof}
Since $Q^{-1}$ is continuous and
$
\mu_t^k\rightharpoonup\mu_t,
$
we have
\[
\mu_O^k\rightharpoonup\mu_O .
\]
Indeed, for every $\varphi\in C(\Omega_O)$,
\[
\int_{\Omega_O}\varphi(y)\,\diff\mu_O^k(y)
=
\int \varphi(Q^{-1}(z))\,\diff\mu_t^k(z)
\longrightarrow
\int \varphi(Q^{-1}(z))\,\diff\mu_t(z)
=
\int_{\Omega_O}\varphi(y)\,\diff\mu_O(y).
\]

Let
$
K_O
$
be a common compact subset of $\Omega_O$.
Consider the far-field transport cost
\[
c(x,y):=-\log(1-x\cdot y),
\qquad
(x,y)\in\Omega_I\times K_O .
\]
Since $\Omega_I$ and $K_O$ are compact and the far-field geometry ensures that
\[
\inf_{x\in\Omega_I,\;y\in K_O}(1-x\cdot y)>0,
\]
the cost
$
c(x,y)
$
is continuous on the compact set $\Omega_I\times K_O$. Hence, it is bounded and uniformly continuous.
Set
\[
u_k:=\log\rho_k,
\qquad
u:=\log\rho .
\]
By the optimal transport characterization of the far-field reflector, $u_k$ is a $c$-concave Kantorovich potential associated with the transport
from $\mu_I$ to $\mu_O^k$. Thus, for some function $\psi_k$ on $K_O$,
\[
u_k(x)
=
\inf_{y\in K_O}
\{c(x,y)-\psi_k(y)\}.
\]
Therefore, for any $x,x'\in\Omega_I$,
\[
|u_k(x)-u_k(x')|
\le
\sup_{y\in K_O}
|c(x,y)-c(x',y)|.
\]
Since $c$ is uniformly continuous on $\Omega_I\times K_O$,
the family $\{u_k\}$ is equicontinuous on $\Omega_I$.

The normalization
$
\rho_k(x_0)=1
$
gives
$
u_k(x_0)=0.
$
Hence
\[
|u_k(x)|
=
|u_k(x)-u_k(x_0)|
\le
\sup_{y\in K_O}
|c(x,y)-c(x_0,y)|,
\]
which shows that $\{u_k\}$ is uniformly bounded.
By the Arzel\`a--Ascoli theorem, every subsequence of $\{u_k\}$ contains a further subsequence converging uniformly on $\Omega_I$. Let such a subsequence, still denoted by $\{u_k\}$, satisfy
\[
u_k\longrightarrow\bar u
\qquad\text{uniformly on }\Omega_I .
\]
Define the $c$-transform
\[
u_k^c(y)
:=
\inf_{x\in\Omega_I}
\{c(x,y)-u_k(x)\},
\qquad y\in K_O.
\]
Then
\[
\|u_k^c-\bar u^c\|_{L^\infty(K_O)}
\le
\|u_k-\bar u\|_{L^\infty(\Omega_I)}
\longrightarrow0 .
\]

Since $(u_k,u_k^c)$ is an optimal Kantorovich dual pair for
$(\mu_I,\mu_O^k)$,
\[
\int_{\Omega_I}u_k\,\diff\mu_I
+
\int_{K_O}u_k^c\,\diff\mu_O^k
=
\mathcal C(\mu_I,\mu_O^k),
\]
where $\mathcal C$ denotes the optimal transport cost associated with $c$.
Because $c$ is bounded and continuous on $\Omega_I\times K_O$, the optimal transport value is continuous under weak convergence of the marginals. Hence
\[
\mathcal C(\mu_I,\mu_O^k)
\longrightarrow
\mathcal C(\mu_I,\mu_O).
\]
Using the uniform convergence of $u_k$ and $u_k^c$, together with
$\mu_O^k\rightharpoonup\mu_O$, we obtain
\[
\int_{\Omega_I}\bar u\,\diff\mu_I
+
\int_{K_O}\bar u^c\,\diff\mu_O
=
\mathcal C(\mu_I,\mu_O).
\]
Thus $\bar u$ is an optimal Kantorovich potential for the limiting
transport problem from $\mu_I$ to $\mu_O$.

By the equivalence between the far-field reflector and the corresponding optimal transport problem, i.e., Theorem \ref{OTRef}
\[
\bar\rho:=e^{\bar u}
\]
is the radial function of a weak reflector associated with $\mu_O$.
By uniqueness of the weak reflector and $\bar\rho(x_0)=1.$, i.e., Theorem \ref{ExiUniThe}, we have
\[
\bar\rho=\rho,
\qquad\text{and hence}\qquad
\bar u=u.
\]
Thus every uniformly convergent subsequence of $\{u_k\}$ has the
same limit $u$. It follows that
\[
u_k\longrightarrow u
\qquad\text{uniformly on }\Omega_I.
\]
Since $\{u_k\}$ is uniformly bounded and the exponential function is uniformly Lipschitz, 
\[
\|\rho_k-\rho\|_{L^\infty(\Omega_I)}
\longrightarrow0.
\]
Hence
\[
\rho_k\longrightarrow\rho
\qquad\text{uniformly on }\Omega_I.
\]
Finally,
\[
R_{\rho_k}
=
\{\rho_k(x)x:\ x\in\Omega_I\},
\qquad
R_\rho
=
\{\rho(x)x:\ x\in\Omega_I\}.
\]
For every $x\in\Omega_I$,
\[
|\rho_k(x)x-\rho(x)x|
=
|\rho_k(x)-\rho(x)|.
\]
Hence
\[
d_H(R_{\rho_k},R_\rho)
\le
\|\rho_k-\rho\|_{L^\infty(\Omega_I)}
\longrightarrow0,
\]
which completes the proof.
\end{proof}
This theorem shows that perturbations of the target measure lead to stable perturbations of the corresponding reflector surfaces.

\begin{lem}[Stability of the reflector-induced measure]\label{lem:reflector_measure_stability}
Let $R_{\rho_k}$ and $R_\rho$ be far-field reflectors such that
\[
\rho_k\to\rho
\qquad\text{uniformly on }\Omega_I.
\]
Let $\mu_{R_k}$ and $\mu_R$ be the corresponding reflector-induced measures given by \eqref{pullbackMea}, i.e.,
\[
\mu_{R_k}(B)
=
\int_{V_k(B)} f(x)\,\diff\sigma(x),
\qquad
\mu_R(B)
=
\int_{V(B)} f(x)\,\diff\sigma(x),
\]
for $B\in \mathcal B(\Omega_O)$.
Then
\[
\mu_{R_k}\rightharpoonup\mu_R
\qquad\text{as }k\to\infty.
\]
\end{lem}

\begin{proof}
Note that the supporting relation is stable under the uniform convergence of the radial functions. Indeed, suppose that
\[
x_k\to x,\qquad y_k\to y,
\qquad x_k\in V_k(y_k).
\]
Then there exists $d_k>0$ such that $p_{y_k,d_k}$ supports $\rho_k$ at
$x_k$. Hence
\[
d_k=\rho_k(x_k)(1-x_k\cdot y_k)
\longrightarrow
\rho(x)(1-x\cdot y)=:d>0.
\]
Moreover,
\[
\rho_k(z)\le p_{y_k,d_k}(z),
\qquad z\in\Omega_I,
\]
with equality at $z=x_k$. Passing to the limit gives
\[
\rho(z)\le p_{y,d}(z),
\qquad z\in\Omega_I,
\]
with equality at $z=x$. Thus $p_{y,d}$ supports $\rho$ at $x$, and hence
$
x\in V(y).
$

Let now $F\subset\Omega_O$ be closed and let $U\subset\Omega_I$ be an
open set containing $V(F)$. We claim that
$
V_k(F)\subset U
$
for all sufficiently large $k$.
Otherwise, there exist a subsequence and points
$x_k\in V_k(F)\setminus U$. By compactness,
\[
x_k\to x\in\Omega_I\setminus U.
\]
For each $k$, choose $y_k\in F$ such that $x_k\in V_k(y_k)$.
Again by compactness,
$
y_k\to y\in F.
$
The stability of the supporting relation proved above gives
\[
x\in V(y)\subset V(F)\subset U,
\]
which is a contradiction.
Thus,
\[
\limsup_{k\to\infty}\mu_{R_k}(F)
\le
\mu_I(U).
\]
Taking the infimum over all open sets $U\supset V(F)$ and using the outer
regularity of $\mu_I$, we obtain
\[
\limsup_{k\to\infty}\mu_{R_k}(F)
\le
\mu_I(V(F))
=
\mu_R(F).
\]
Finally,
$
\mu_{R_k}(\Omega_O)
=
\mu_R(\Omega_O)
=
\mu_I(\Omega_I).
$
Hence, by the Portmanteau theorem,
$
\mu_{R_k}\rightharpoonup\mu_R.
$
\end{proof}

\begin{cor}\label{cor:GOA_measure_stability}
Let $\mu_t^k\rightharpoonup\mu_t$ and let $R_{\rho_k}$ and $R_\rho$ be the corresponding far-field reflectors satisfying the assumptions of Theorem~\ref{StaTarDomThe} and Lemma~\ref{lem:reflector_measure_stability}.
Define
\[
\mu_{\mathrm{GOAS}}^k
:=
Q_\sharp\mu_{R_k},
\qquad
\mu_{\mathrm{GOAS}}
:=
Q_\sharp\mu_R.
\]
Then
\[
\mu_{\mathrm{GOAS}}^k
\rightharpoonup
\mu_{\mathrm{GOAS}}
\qquad\text{as }k\to\infty.
\]
\end{cor}

\begin{proof}
By Theorem~\ref{StaTarDomThe} and Lemma~\ref{lem:reflector_measure_stability}, we have
\[
\mu_{R_k}\rightharpoonup\mu_R.
\]
Since $Q$ is continuous, weak convergence is preserved under push-forward, and hence
$
Q_\sharp\mu_{R_k}
\rightharpoonup
Q_\sharp\mu_R.
$
\end{proof}

\begin{rem}[Well-posedness of GOAS measure]
The preceding results establish the well-posedness of the GOAS measure. Indeed, for prescribed source and target measures, the far-field reflector exists and is unique up to a positive scaling from Theorem \ref{ExiUniThe}. Therefore, $\mu_{\mathrm{GOAS}}$ is uniquely defined. Moreover, Corollary~\ref{cor:GOA_measure_stability} shows that $\mu_{\mathrm{GOAS}}$ depends continuously, in the sense of weak convergence, on perturbations of the target measure.

This stability is particularly relevant to the numerical realization developed in the next section, where the supporting hyperparaboloid construction requires a discrete approximation of the target measure. Hence different consistent discretization strategies, including density-based and sample-based constructions, lead to stable GOAS measures.
\end{rem}

\section{Implementation of far-field GOAS}
\label{sec:GOASAlg}

The previous section establishes the reflector-induced transport underlying far-field GOAS. We now describe one numerical realization of this framework. The implementation consists of an offline stage, in which a reflector is constructed, and an online stage, in which independent source samples are transported through the constructed reflector.

In this work, the offline reflector problem is solved by a supporting hyperparaboloid method. This is only one possible numerical realization of far-field GOAS; the sampling framework itself is not tied to a particular reflector solver. Once the reflector has been constructed, we consider two sampling mechanisms. The first uses the primal--dual reflector geometry directly for the nonsmooth piecewise reflector. The second applies a softmin regularization and uses the resulting smooth reflection map.

\subsection{Supporting hyperparaboloid construction}
\label{sec:SPmethod}
Let
$
E:=\mu_I(\Omega_I)
$
denote the total source energy.
The supporting hyperparaboloid construction requires a finite weighted
approximation of the target measure,
\begin{align}
\label{TarDirDis}
\mu_t^K
=
\sum_{i=1}^K w_i\delta_{z_i},
\qquad
w_i>0,
\qquad
\sum_{i=1}^K w_i=E.
\end{align}
The particular construction of \eqref{TarDirDis} is not essential to the
reflector solver. If a target density is available, the points and weights may
be obtained from quadrature, random or low-discrepancy points, or importance
sampling. If the target is available through samples, one may use the
empirical measure directly or compress a large sample set into weighted
representative points, for example by clustering. Thus, the numerical
construction below only requires the weighted pairs
$\{(z_i,w_i)\}_{i=1}^K$.

The discrete target measure is transferred to the output sphere through
$Q^{-1}$:
\begin{align}
\label{TarDirDisSph}
\mu_O^K
=
(Q^{-1})_\sharp\mu_t^K
=
\sum_{i=1}^K w_i\delta_{y_i},
\qquad
y_i=Q^{-1}(z_i).
\end{align}
Since \eqref{TarDirDis} is a discrete measure, its masses $w_i$ are unchanged
under the push-forward.

For the output directions $\{y_i\}_{i=1}^K$, Definition \ref{def:SP} and Proposition \ref{rhoTconvexRef} imply that a piecewise-hyperparaboloidal reflector can be represented as
$$
R_{\tilde\rho}(d):=\{x\tilde\rho(x):x\in\Omega_I\},
$$
where
\begin{align}
\label{eq:nonsmoothreflector}
\tilde\rho(x)
=
\min_{1\le i\le K}
\Phi_i(x),
\qquad
\Phi_i(x)
:=
\frac{d_i}{1-x\cdot y_i},
\end{align}
where  \(d=(d_1,\ldots,d_K)\in(0,\infty)^K\) is the vector of focal parameters of the supporting hyperparaboloids.
The corresponding visibility cell is
\begin{align}
\label{eq:visibility_cell}
V_i(d)
:=
\left\{
x\in\Omega_I:
\Phi_i(x)
=
\min_{1\le j\le K}\Phi_j(x)
\right\}.
\end{align}
The source mass reflected toward $y_i$ is therefore
\begin{align}
\label{eq:cell_mass}
G_i(d)
:=
\mu_I(V_i(d))
=
\int_{V_i(d)}f(x)\,\diff\sigma(x).
\end{align}
Therefore, the discrete measure induced by \(R_{\tilde\rho}(d)\) is
\begin{align}
\label{eq:discrete_measure_bR}
\mu_{R_{\tilde\rho}(d)}
=
\sum_{i=1}^K G_i(d)\,\delta_{y_i}.
\end{align}
Since the visibility cells partition \(\Omega_I\), the induced masses automatically satisfy
\begin{align}\label{eq:induced_mass_conservation}
\sum_{i=1}^K G_i(d)
=
\mu_I(\Omega_I)
=
E.
\end{align}
The discrete reflector problem is to determine
$d=(d_1,\ldots,d_K)$ such that
\begin{align}
\label{eq:discrete_reflector_problem}
G_i(d)=w_i,
\qquad i=1,\ldots,K.
\end{align}

Multiplying all focal parameters by the same positive constant does not change the visibility cells. We therefore fix one focal parameter, for example $d_1$, and determine the remaining parameters. The supporting hyperparaboloid method iteratively adjusts these parameters using the monotonic dependence of the reflector-induced masses $G_i(d)$ on $d_i$ \cite{caffarelli1999,kochengin2003}. The iteration is terminated when
\begin{align}
\label{SPconErr}
 \mathcal E(d):=
 \frac{1}{E\sqrt K}
\left(
\sum_{i=1}^K
\bigl(
 G_i(d)-w_i
\bigr)^2
\right)^{1/2}
\le \epsilon,
\end{align}
where $\epsilon>0$ is the prescribed reflector-construction tolerance.
To retain the update used in our implementation, define
\begin{align}\label{AlgDicon}
    G_i(d)
    \le
    w_i+\frac{\epsilon}{\sqrt{K(K-1)}},
    \qquad i=2,\dots,K,
\end{align}
and fix \(d_1\) to remove the scale invariance. An initial focal vector is chosen as
\begin{align}\label{IniRef}
    d^0=(d_1,\alpha d_2,\dots,\alpha d_K),
    \qquad \alpha>1.
\end{align}
Algorithm~\ref{AlgSP} summarizes the resulting construction. Importantly, the construction uses only the discrete target measure and evaluations of the induced cell masses; no gradient of the target density is required.

\begin{algorithm}[htbp]
\renewcommand{\algorithmicrequire}{\textbf{Input:}}
\renewcommand{\algorithmicensure}{\textbf{Output:}}
\caption{Supporting hyperparaboloid construction}
\label{AlgSP}
\begin{algorithmic}[1]
\Require Dirac output measure \(\mu_O^K\) given by \eqref{TarDirDisSph}, source measure \(\mu_I\), tolerance \(\epsilon>0\)
\Ensure Focal vector \(d=(d_1,\dots,d_K)\) satisfying \eqref{SPconErr}
\State Choose an initial focal vector \(d^0\) as in \eqref{IniRef}; keep \(d_1\) fixed
\State Set \(\widetilde d=d^0\), \(j=0\), and initialize increments \(\Delta d_1=0\) and \(\Delta d_i\), \(i=2,\dots,K\)
\State Evaluate \(G_i(\widetilde d)\), \(i=1,\dots,K\)
\While{\(\mathcal E(\widetilde d)>\epsilon\)}
    \State Let \(J\subset\{2,\dots,K\}\) contain the indices violating \eqref{AlgDicon}
    \If{\(J=\emptyset\)}
        \State Set \(d^{j+1}=\widetilde d\), and $\Delta d_i= d_i^j/3,\;i\in\{2,3, \dots,K\}$
    \Else
        \State Keep \(d^{j+1}=d^j\) and \(\Delta d_i=\Delta d_i/2\) for \(i\in J\)
    \EndIf
    \State Put \(\widetilde d_i=d^{j+1}_i-\Delta d_i,\;i\in\{2,3, \dots,K\}\)
    \State Recompute \(G_i(\widetilde d)\) and set \(j\leftarrow j+1\)
\EndWhile
\State \Return \(d=\widetilde d\)
\end{algorithmic}
\end{algorithm}


\subsubsection{Monte Carlo evaluation of reflector-induced masses}
\label{sec:MC_mass}

For a given focal-parameter vector $d$, the mass reflected toward the output direction
$y_i$ is
\begin{equation}
\label{eq:Gi_farfield}
G_i(d)
=
\int_{\Omega_I}
f(x)\mathbf{1}_{V_i(d)}(x)\,\diff\sigma(x),
\end{equation}
where $V_i(d)$ is the visibility cell associated with the $i$th supporting hyperparaboloid and $\mathbf{1}_{V_i(d)}$ denotes its indicator function. Since an explicit evaluation of these visibility cells is generally
computationally expensive, we approximate \eqref{eq:Gi_farfield} by Monte Carlo ray tracing \cite{shirley2008realistic,fournier2010freeform,sun2024}.

Let
$
X_1,\ldots,X_N
$
be independent source rays drawn from the surface measure $\sigma$ on $\Omega_I$.
For each $X_\ell$, the active supporting hyperparaboloid is determined by
\begin{equation}
\label{IdeSP_farfield}
\iota_d(X_\ell)
:=
\operatorname*{arg\,min}_{1\le j\le K}
\frac{d_j}{1-X_\ell\cdot y_j}.
\end{equation}
The reflector-induced mass $G_i(d)$ is then approximated by
\begin{equation}
\label{ApprG_farfield}
\widetilde G_i(d)
:=
\frac{\sigma(\Omega_I)}{N}
\sum_{\ell=1}^N
f(X_\ell)
\mathbf{1}_{\{\iota_d(X_\ell)=i\}},
\qquad
i=1,\ldots,K.
\end{equation}
If we set $\mu_I$ as uniform source measure on $\Omega_I$, then $f=E/\sigma(\Omega_I)$.

In Algorithm~\ref{AlgSP},
$
G(d)
$
is replaced by its Monte Carlo approximation
\[
\widetilde G(d)
=
\bigl(\widetilde G_1(d),\ldots,\widetilde G_K(d)\bigr).
\]
The same set of source rays is reused throughout the reflector iteration,
so that changes in $\widetilde G(d)$ are caused only by updates of the
focal parameters. Each evaluation of $\widetilde G(d)$ requires
$\mathcal{O}(KN)$ operations.

\subsection{Sampling from the constructed reflector}
\label{sec:sampling_reflector}

Once reflector $\tilde\rho$ has been constructed, approximate target samples can be generated by transporting source samples through the reflector system. Since \eqref{eq:nonsmoothreflector} is continuous but generally nonsmooth across the interfaces of the visibility cells, we consider two complementary sampling realizations.

\subsubsection{Primal--dual consistency resampling}
\label{sec:primal_dual_sampling}

Away from cell interfaces, the active supporting hyperparaboloid is unique.
Define the active-cell index by
$
\iota_d(x)
$
given by \eqref{IdeSP_farfield}.
The discrete reflector assigns a source direction $x$ to the output direction
\begin{align}
\label{dis:refdir}
T_{\widetilde\rho}(x)
=
y_{\iota_d(x)}.
\end{align}
Direct use of \eqref{dis:refdir} produces samples on the discrete support $\{y_i\}_{i=1}^K$. To obtain a continuous numerical realization without differentiating the nonsmooth reflector, we exploit the primal--dual structure developed in Section~\ref{sec:dual_transport}.

For each nonempty visibility cell $V_i$, choose a representative direction $x_i\in V_i$ from the source rays assigned to that cell. Since $y_i$ is active at $x_i$,
\[
\tilde\rho(x_i)
=
\frac{d_i}{1-x_i\cdot y_i}.
\]
We therefore define
\begin{align}
\label{eq:dual_parameter}
d_i^\ast
:=
\frac{1}{\widetilde\rho(x_i)}
=
\frac{1-x_i\cdot y_i}{d_i}.
\end{align}
Using the representative pairs $\{(x_i,y_i)\}_{i=1}^K$, we construct the discrete dual approximation reflector
\begin{align}
\label{eq:dual_discrete_reflector}
\tilde\rho^\ast(y)
=
\min_{1\le i\le K}
\frac{d_i^\ast}{1-x_i\cdot y},
\qquad
y\in\Omega_O,
\end{align}
and define its active-cell index by
\begin{align}
\label{eq:active_dual_index}
\iota_{d^\ast}^\ast(y)
:=
\operatorname*{arg\,min}_{1\le i\le K}
\frac{d_i^\ast}{1-x_i\cdot y}.
\end{align}

The primal index identifies the output cell associated with a source
direction, whereas the dual index provides a reciprocal consistency check.
Given a source sample $x'$, let
\[
j=\iota_d(x').
\]
A candidate $z$ is generated locally around $z_j=Q(y_j)$ and mapped to
$\widetilde y=Q^{-1}(z)$. The candidate is retained whenever
\[
\iota_{d^\ast}^\ast(\widetilde y)=j.
\]
Repeating this procedure gives a continuous, primal--dual consistent numerical realization of the discrete reflector transport, see Algorithm \ref{alg:dualresim}.

\begin{algorithm}[htbp]
\renewcommand{\algorithmicrequire}{\textbf{Input:}}
\renewcommand{\algorithmicensure}{\textbf{Output:}}
\caption{Primal--dual consistency resampling}
\label{alg:dualresim}
\begin{algorithmic}[1]

\Require
$\tilde\rho$, $\tilde\rho^\ast$, source measure $\mu_I$,
transformation $Q$, number of samples $M$

\Ensure
Approximate target samples $\{z_m'\}_{m=1}^M$

\For{$m=1,\ldots,M$}
    \State Draw $x_m'\sim\mu_I$ and set
    $j=\iota_d(x_m')$.
    \Repeat
        \State Draw a local candidate $z$ around $Q(y_j)$.
        \State Set $\widetilde y=Q^{-1}(z)$.
    \Until{$\iota_{d^\ast}^\ast(\widetilde y)=j$}
    \State Set $z_m'=z$.
\EndFor

\end{algorithmic}
\end{algorithm}

In Algorithm \ref{alg:dualresim}, the local candidate may be generated from a simple proposal centered at $Q(y_j)$, such as a Gaussian or a uniform distribution. The particular choice of the local proposal is not essential to the primal--dual consistency mechanism.

The resulting approximate target samples preserve the primal--dual cell correspondence, while their local within-cell distribution is determined by the chosen proposal mechanism.
The local proposal explores the target associated with the selected primal cell, while the dual selector enforces reciprocal cell consistency. The procedure avoids the evaluation or interpolation of surface normals at the nonsmooth interfaces. 

\subsubsection{Softmin smoothing and direct reflection maps}
\label{sec:softmin_maps}

The second sampling realization replaces the piecewise reflector by a smooth approximation. For the discrete reflector $\tilde\rho$ given by \eqref{eq:nonsmoothreflector}, define
\begin{align}
\label{eq:rhodisSmooth}
\tilde\rho_\lambda(x)
=
-\lambda
\log\left(
\sum_{i=1}^K
\exp\left[-\frac{\Phi_i(x)}{\lambda}\right]
\right),
\end{align}
The parameter $0<\lambda<1$ controls the balance between smoothness and approximation accuracy.

\begin{thm}[Softmin smoothing]\label{thm:rhoSmo}
Fix $x\in\Omega_I$ and let
\[
I^\ast(x)
=
\{i:\ \Phi_i(x)=\tilde\rho(x)\}
\]
be the set of active minimizers.
Let $m=|I^\ast(x)|\ge 1$ and assume that there exists $\tau>0$ such that
\[
\Phi_j(x)-\tilde\rho(x)\ge \tau,
\qquad j\notin I^\ast(x).
\]
Then
\begin{align}\label{ineq:rhoappro}
0
\le
\tilde\rho(x)-\tilde\rho_\lambda(x)
\le
\lambda\log m
+
C_1\lambda e^{-\tau/\lambda},
\end{align}
where $C_1>0$ depends only on $K$ and $m$.
Moreover, if the minimizer is unique, i.e., $m=1$, then $\tilde\rho$ is differentiable at $x$ and
\begin{align}\label{ineq:Drhoappro}
\left|D\tilde\rho_\lambda(x)-D\tilde\rho(x)\right|
\le
C_2 e^{-\tau/\lambda},
\end{align}
where $C_2>0$ depends on the functions $\{\Phi_i\}_{i=1}^K$.
Here $D$ denotes differentiation with respect to $x$ in a local coordinate chart on $\Omega_I$.
\end{thm}

\begin{proof}
Let
\[
\tilde\rho(x)=\min_{1\le i\le K}\Phi_i(x).
\]
For $i\in I^\ast(x)$, we have $\Phi_i(x)=\tilde\rho(x)$.
For $j\notin I^\ast(x)$, write
\[
\tau_j=\Phi_j(x)-\tilde\rho(x)\ge \tau>0.
\]
Then
\begin{align*}
\tilde\rho_\lambda(x)
&=
-\lambda\log
\left[
\sum_{i=1}^{K}
e^{-\Phi_i(x)/\lambda}
\right]  \\
&=
-\lambda\log
\left[
m e^{-\tilde\rho(x)/\lambda}
+
\sum_{j\notin I^\ast(x)}
e^{-\Phi_j(x)/\lambda}
\right]  \\
&=
\tilde\rho(x)
-
\lambda
\log
\left[
m+
\sum_{j\notin I^\ast(x)}
e^{-\tau_j/\lambda}
\right].
\end{align*}
Hence
\[
\tilde\rho(x)-\tilde\rho_\lambda(x)
=
\lambda
\log
\left[
m+
\sum_{j\notin I^\ast(x)}
e^{-\tau_j/\lambda}
\right].
\]
Using $\tau_j\ge\tau$, we obtain
\begin{align*}
\tilde\rho(x)-\tilde\rho_\lambda(x)
&\le
\lambda
\log
\left[
m+(K-m)e^{-\tau/\lambda}
\right]  \\
&=
\lambda\log m
+
\lambda
\log
\left[
1+\frac{K-m}{m}e^{-\tau/\lambda}
\right]  \\
&\le
\lambda\log m
+
\frac{\lambda(K-m)}{m}e^{-\tau/\lambda}.
\end{align*}
This proves \eqref{ineq:rhoappro}.

If the minimizer is unique, let $I^\ast(x)=\{i^\ast\}$.
Then
\[
\tilde\rho(x)=\Phi_{i^\ast}(x),
\qquad
D\tilde\rho(x)=D\Phi_{i^\ast}(x).
\]
Moreover,
\[
D\tilde\rho_\lambda(x)
=
\frac{
\sum_{i=1}^{K}
e^{-\Phi_i(x)/\lambda}D\Phi_i(x)
}{
\sum_{i=1}^{K}
e^{-\Phi_i(x)/\lambda}
}.
\]
Therefore,
\begin{align*}
D\tilde\rho_\lambda(x)-D\tilde\rho(x)
&=
\frac{
\sum_{i\neq i^\ast}
e^{-\Phi_i(x)/\lambda}
\left(D\Phi_i(x)-D\Phi_{i^\ast}(x)\right)
}{
\sum_{i=1}^{K}
e^{-\Phi_i(x)/\lambda}
}.
\end{align*}
Since
\[
e^{-\Phi_i(x)/\lambda}
=
e^{-\Phi_{i^\ast}(x)/\lambda}
e^{-(\Phi_i(x)-\Phi_{i^\ast}(x))/\lambda}
\le
e^{-\Phi_{i^\ast}(x)/\lambda}e^{-\tau/\lambda},
\qquad i\neq i^\ast,
\]
we obtain
\begin{align*}
\left|D\tilde\rho_\lambda(x)-D\tilde\rho(x)\right|
&\le
e^{-\tau/\lambda}
\sum_{i\neq i^\ast}
\left|D\Phi_i(x)-D\Phi_{i^\ast}(x)\right|  \\
&\le
C_2 e^{-\tau/\lambda}.
\end{align*}
This proves \eqref{ineq:Drhoappro}.
\end{proof}

The smooth reflector induces a single-valued reflection map.
We use a local coordinate system on $S^n$, namely
\begin{align}\label{sys:coor}
x=(t^1,\dots,t^n,s),\qquad 
s=\sqrt{1-|t|^2},\qquad 
t=(t^1,\dots,t^n)\in\mathbb R^n,\quad |t|<1.
\end{align}

\begin{prop}[Approximate reflector map]\label{prop:softmin_reflector_map}
The softmin-smoothed reflector $R_{\tilde\rho_\lambda}$ induces a single-valued reflector map
\[
T_{\tilde\rho_\lambda}:\Omega_I\to \Omega_O
\]
given by
\begin{align}\label{eq:softmin_reflection_map}
T_{\tilde\rho_\lambda}(x)
=
x+
\frac{
2\tilde\rho_\lambda
\bigl(
D\tilde\rho_\lambda
-
x(\tilde\rho_\lambda+D\tilde\rho_\lambda\cdot x)
\bigr)
}{
\tilde\rho_\lambda^2
+
|D\tilde\rho_\lambda|^2
-
(D\tilde\rho_\lambda\cdot x)^2
}.
\end{align}
Here $D$ denotes differentiation with respect to $x$ in \eqref{sys:coor}, and
\begin{align}\label{eq:softmin_Drho}
D\tilde\rho_\lambda(x)
=
\frac{
\sum_{i=1}^{K}
\exp\!\left(-\frac{\Phi_i(x)}{\lambda}\right)
D\Phi_i(x)
}{
\sum_{i=1}^{K}
\exp\!\left(-\frac{\Phi_i(x)}{\lambda}\right)
},
\end{align}
where,
\begin{align}\label{eq:Df}
D\Phi_i(x)
=
\frac{d_i}{(1-x\cdot y_i)^2}
\left(
y_i-\frac{y_{i,n+1}}{x_{n+1}}x
\right).
\end{align}
\end{prop}

\begin{proof}
The smooth reflector is represented as the radial graph
\[
R_{\tilde\rho_\lambda}
=
\{x\tilde\rho_\lambda(x):x\in\Omega_I\}.
\]
The outward unit normal vector of a radial graph over the sphere can be written as
\begin{align}\label{eq:unitnormal_softmin}
\nu_\lambda
=
\frac{
D\tilde\rho_\lambda
-
x(\tilde\rho_\lambda+D\tilde\rho_\lambda\cdot x)
}{
\sqrt{
\tilde\rho_\lambda^2
+
|D\tilde\rho_\lambda|^2
-
(D\tilde\rho_\lambda\cdot x)^2
}
}.
\end{align}
Taking the inner product with $x$ gives
\[
x\cdot\nu_\lambda
=
-\frac{\tilde\rho_\lambda}{
\sqrt{
\tilde\rho_\lambda^2
+
|D\tilde\rho_\lambda|^2
-
(D\tilde\rho_\lambda\cdot x)^2
}
}.
\]
By the physical law of reflection, the reflected direction is
\[
T_{\tilde\rho_\lambda}(x)
=
x-2(x\cdot\nu_\lambda)\nu_\lambda.
\]
Substituting \eqref{eq:unitnormal_softmin} into this identity yields \eqref{eq:softmin_reflection_map}.
In the local chart \eqref{sys:coor}, we have
\[
D\Phi_i(x)
=
\frac{d_i}{(1-x\cdot y_i)^2}
\left(
y_i-\frac{y_{i,n+1}}{x_{n+1}}x
\right).
\]
\end{proof}

The far-field duality also associates a dual reflector with
$\tilde\rho_\lambda$:
\begin{align}
\label{eq:smooth_dual}
\widetilde\rho_\lambda^\ast(y)
=
\inf_{x\in\Omega_I}
\frac{1}
{\tilde\rho_\lambda(x)(1-x\cdot y)}.
\end{align}
Whenever the reflection maps are single-valued and one-to-one, the reciprocal relation gives
\begin{align}
T_\lambda^\ast
=
T_\lambda^{-1}.
\end{align}
The formula \eqref{eq:softmin_reflection_map} follows from the unit normal of the radial graph and the law of reflection. The inverse relation is an immediate consequence of Proposition~\ref{LTDual} applied to the smooth pair $(\tilde\rho_\lambda,\tilde\rho_\lambda^\ast)$.
Thus the smooth formulation provides an explicit forward transport map and, through the far-field duality, its natural inverse map.


\subsection{The far-field GOAS algorithm}\label{sec:GOASalgorithm}

The complete workflow is summarized in Algorithm~\ref{alg:GOAS}. In this work, the supporting hyperparaboloid method is adopted as the numerical solver for constructing the far-field reflector; other solvers for the same reflector problem may be incorporated without changing the overall GOAS framework.
We use this construction because it requires no gradient information of the target density and is applicable in general dimensions.

\begin{algorithm}[htbp]
\renewcommand{\algorithmicrequire}{\textbf{Input:}}
\renewcommand{\algorithmicensure}{\textbf{Output:}}
\caption{Far-field Geometric Optics Approximation Sampling (GOAS)}
\label{alg:GOAS}
\begin{algorithmic}[1]

\Require
Target information (density, samples, or weighted particles),
source measure $\mu_I$, transformation $Q$, reflector tolerance $\epsilon$,
number of samples $M$, and optionally a smoothing parameter $\lambda$

\Ensure
Approximate target samples $\{z_m'\}_{m=1}^M$

\State \textbf{Offline stage:}
Construct a discrete approximation $\mu_t^K$ using \eqref{TarDirDis} required by the chosen reflector solver

\State Construct a far-field discrete reflector; in this work use supporting hyperparaboloid (Algorithm~\ref{AlgSP}) to obtain $\tilde\rho$

\State \textbf{Online sampling stage:}

\State Choose one of the following realizations:

\State \hspace{0.5cm}
\textbf{(a) Nonsmooth realization:}
apply Algorithm~\ref{alg:dualresim}.

\State \hspace{0.5cm}
\textbf{(b) Smooth realization:}
construct $\tilde\rho_\lambda$ from \eqref{eq:rhodisSmooth},
draw $x_m\sim\mu_I$, and set
$$
z_m'=Q(T_{\tilde\rho_\lambda}(x_m))\;,m=1,\ldots,M.
$$

\State \Return $\{z_m'\}_{m=1}^M$.

\end{algorithmic}
\end{algorithm}

The two realizations share the same offline reflector construction. The primal--dual consistency resampling strategy operates directly on the nonsmooth piecewise reflector, whereas the softmin realization yields an explicit smooth transport map. Neither the construction of the discrete reflector nor the subsequent sampling procedure requires gradients of the target density. In the next section, we analyze the principal approximation errors arising from target discretization, Monte Carlo evaluation of the reflector-induced masses, the reflector construction tolerance, and softmin smoothing.




\section{Wasserstein error analysis} \label{sec:wasserstein_error}

In this section, we quantify the error of the far-field GOAS measure in Wasserstein distance.  The analysis separates three errors generated during the offline construction and the smooth-map realization: (i) approximation of the target measure by a weighted discrete measure, (ii) the reflector-construction error caused by the stopping tolerance and the Monte Carlo evaluation of the reflector-induced masses, and (iii) the softmin smoothing error.
The additional within-cell proposal used in the primal–dual consistency resampling procedure in Section \ref{sec:primal_dual_sampling} is not included in the present analysis.

\begin{defn}[Wasserstein distances \cite{villani2009,villani2021,benamou2020}]
Given two compact manifold $\mathcal{X}$ and $\mathcal{Y}$  endowed with a continuous, bounded from below cost function $c(x,y):\mathcal{X}\times \mathcal{Y}\to\mathbb{R}$ for transporting one unit of mass from $x$ to $y$.
Let $\mu_1$ and $\mu_2$ be measures on $\mathcal{X}$, $\mathcal{Y}$ respectively, the optimal transportation cost
between $\mu_1$ and $\mu_2$ is defined by
\begin{align}\label{WasDis}
OT_c(\mu_1,\mu_2)=\inf_{\gamma\in\Gamma(\mu_1,\mu_2)}\int_{\mathcal{X}\times \mathcal{Y}}c(x,y)\diff\gamma(x,y)
\end{align}
where $\Gamma(\mu_1,\mu_2)$ is the set of measures $\gamma$ on $\mathcal{X}\times \mathcal{Y}$, called couplings of $\mu_1$ and $\mu_2$, satisfying $\gamma(A\times \mathcal{Y})=\mu_1(A)$ and $\gamma(\mathcal{X}\times B)=\mu_2(B)$ for all Borel subsets $A\subset \mathcal{X}$ and $B\subset\mathcal{Y}$. In particular, if $\mathcal{X}=\mathcal{Y}$ and $(\mathcal{X},d_{\mathcal{X}})$ is a compact metric space, then for $p\in [1,\infty]$, the \emph{Wasserstein distance} of order $p$ between $\mu_1$ and $\mu_2$ is defined as $W_p(\mu_1,\mu_2)=(OT_c(\mu_1,\mu_2))^{\frac{1}{p}}$ with
$$
c(x,y)=d^p_{\mathcal{X}}(x,y),\quad x,y\in\mathcal{X}.
$$
\end{defn}
\begin{rem}
Given a measure $\gamma\in\Gamma(\mu_1,\mu_2)$ and pair of location $(x,y)$, the value of $\gamma(x,y)$ tells us what proportion of mass of $\mu_1$ at $x$ ought to be transferred to $y$ in order to reconfigure $\mu_1$ into $\mu_2$.
In this section, the metric on the sphere is taken as $d_{S^n}$.
\end{rem}

\subsection{Target-measure discretization}\label{sec:target_Wp_error}

Let
$
\mu_t^K
$
be the discrete target measure introduced in \eqref{TarDirDis}, and define
\begin{align}
\label{eq:target_dis_error}
\delta_K
:=
W_p(\mu_t^K,\mu_t).
\end{align}
The quantity $\delta_K$ measures only the accuracy of the target representation and is independent of the reflector solver.

\begin{prop}[Density-based target discretization]
\label{prop:density_target_discretization}
Let $\{z_i,a_i\}_{i=1}^K$ be a quadrature or design rule, and let
\begin{align}\label{eq:density_based_weights}
\mu_t^K=\sum_{i=1}^K w_i\delta_{z_i},
\qquad
w_i=E
\frac{\widetilde w_i}
{\sum_{j=1}^K \widetilde w_j},
\qquad
\widetilde w_i=a_i\pi(z_i)
\end{align}
be the density-based approximation of $\mu_t$ given by \eqref{TarMea}.
Then there exists a measurable partition $\{A_i\}_{i=1}^K$ of $\Omega$
such that
\[
\delta_K
\le
\left(
\sum_{i=1}^K
\int_{A_i}|z-z_i|^p\,\diff\mu_t(z)
\right)^{1/p}.
\]
\end{prop}

\begin{proof}
Since $\mu_t$ is absolutely continuous, it is non-atomic. Hence there
exists a measurable partition $\{A_i\}_{i=1}^K$ of $\Omega$ satisfying
\[
\mu_t(A_i)=w_i,
\qquad i=1,\ldots,K.
\]
Define $F_K(z)=z_i$ for $z\in A_i$. Then
\[
(F_K)_\sharp\mu_t
=
\sum_{i=1}^K w_i\delta_{z_i}
=
\mu_t^K.
\]
Thus $(\operatorname{Id},F_K)_\sharp\mu_t$ is an admissible coupling,
and therefore
\[
W_p^p(\mu_t,\mu_t^K)
\le
\int_\Omega |z-F_K(z)|^p\,\diff\mu_t(z)
=
\sum_{i=1}^K
\int_{A_i}|z-z_i|^p\,\diff\mu_t(z).
\]
\end{proof}
The same estimate applies to the importance-weighted construction $\widetilde w_i=\pi(z_i)/q(z_i)$ where \(q\) is proposal density.

\begin{prop}[Sample-based target discretization]
\label{prop:sample_target_discretization}
Let $z_1,\ldots,z_K$ be samples from $\mu_t$. Define
\[
\mu_t^K
=
\frac EK\sum_{i=1}^K\delta_{z_i}.
\]
Then
\[
\delta_K
=
W_p(\mu_t^K,\mu_t)
\longrightarrow0
\qquad\text{almost surely as }K\to\infty.
\]
\end{prop}

\begin{proof}
The empirical measures satisfy
$\mu_t^K\rightharpoonup\mu_t$ almost surely. Moreover, by the strong
law of large numbers,
\[
\int_\Omega |z|^p\,\diff\mu_t^K(z)
=
\frac1K\sum_{i=1}^K|z_i|^p
\longrightarrow
\int_\Omega |z|^p\,\diff\mu_t(z)
\]
almost surely.
The characterization of $p$-Wasserstein convergence therefore gives the result.
\end{proof}

Moreover, if $\mu_t$ has sufficiently high finite moments, then \cite{fournier2015} gives the non-asymptotic estimate
\[
\mathbb E\!\left[W_p^p(\mu_t^K,\mu_t)\right]
\le C
\begin{cases}
K^{-1/2}, & n<2p,\\[1mm]
K^{-1/2}\log(1+K), & n=2p,\\[1mm]
K^{-p/n}, & n>2p,
\end{cases}
\]
where $C>0$ is independent of $K$ but may depend on $n$, $p$, and the moments of $\mu_t$.

%

\subsection{Reflector construction error}\label{sec:reflector_mass_error}

\begin{lem}\label{lem:mass_perturbation}
Let
\[
\nu_1=\sum_{i=1}^K a_i\delta_{y_i},
\qquad
\nu_2=\sum_{i=1}^K b_i\delta_{y_i}
\]
be measures supported on the same points in $\Omega_O$, and let
$
D_O:=\operatorname{diam}_{S^n}(\Omega_O).
$
Then
\begin{align}
\label{eq:mass_perturbation}
W_p(\nu_1,\nu_2)
\le
D_O
\left(
\frac12
\sum_{i=1}^K|a_i-b_i|
\right)^{1/p}.
\end{align}
\end{lem}

\begin{proof}
Put \(c_i:=\min\{a_i,b_i\}\),
\(r_i:=a_i-c_i\), and \(s_i:=b_i-c_i\).  Then
\[
    m:=\sum_{i=1}^K r_i
      =\sum_{i=1}^K s_i
      =\frac12\sum_{i=1}^K|a_i-b_i|.
\]
If \(m=0\), the result is immediate.  If \(m>0\), define
\[
\gamma
:=
\sum_{i=1}^K c_i\delta_{(y_i,y_i)}
+
\frac1m
\sum_{i,j=1}^K
r_i s_j\,\delta_{(y_i,y_j)}.
\]
The first and second marginals of \(\gamma\) are \(\nu_1\) and \(\nu_2\), respectively.  The common mass \(c_i\) is transported at zero cost, whereas the remaining mass \(m\) travels a distance at most \(D_O\).  Thus,
\[
    W_p(\nu_1,\nu_2)^p
    \leq\int_{\Omega_O\times\Omega_O}d^p_{S^n}(y,y’)\,\mathrm d\gamma(y,y’)
    \leq D_O^p m,
\]
which proves~\eqref{eq:mass_perturbation}.
\end{proof}

In practice, $G_i(d)$ is approximated by the Monte Carlo estimator introduced in Section~\ref{sec:MC_mass},
\[
\widetilde G_i(d)
=
\frac{E}{N}
\sum_{\ell=1}^N
\mathbf 1_{\{\iota_d(X_\ell)=i\}},
\qquad
X_\ell\sim\frac{\mu_I}{E}.
\]
In Algorithm \ref{AlgSP} and \eqref{SPconErr}, the actual stopping criterion is
\begin{equation}
\frac{1}{\sqrt K}
\left(
\sum_{i=1}^K
\bigl(
\widetilde G_i(d)-w_i
\bigr)^2
\right)^{1/2}
\le E\epsilon.
\label{eq:computed_stopping_criterion}
\end{equation}
For a fixed focal vector $d$, define
\[
\eta_N(d)
:=
\sum_{i=1}^K
\left|
\widetilde G_i(d)-G_i(d)
\right|.
\]

\begin{lem}\label{lem:MC_mass_error}
For every fixed $d$,
\begin{align}
\label{eq:MC_L1_bound}
\mathbb E\,\eta_N(d)
\le
E\sqrt{\frac{K}{N}}.
\end{align}
\end{lem}

\begin{proof}
For each $i$, define
\[
p_i
:=
\frac{G_i(d)}{E}.
\]
Since $X_\ell\sim\mu_I/E$,
\[
\mathbb P\bigl(\iota_d(X_\ell)=i\bigr)
=
p_i.
\]
Hence
$
\mathbf 1_{\{\iota_d(X_\ell)=i\}}
$
is a Bernoulli random variable with parameter $p_i$, and therefore
\[
\mathbb E\bigl[\widetilde G_i(d)\bigr]
=
G_i(d).
\]
Moreover,
\begin{align*}
\mathbb E
\left[
\bigl(
\widetilde G_i(d)-G_i(d)
\bigr)^2
\right]
=
\frac{E^2}{N}p_i(1-p_i)
=
\frac{G_i(d)\bigl(E-G_i(d)\bigr)}{N}.
\end{align*}
By the Cauchy--Schwarz inequality,
\begin{align*}
\mathbb E\,\eta_N(d)
&=
\sum_{i=1}^K
\mathbb E
\left|
\widetilde G_i(d)-G_i(d)
\right|
\\
&\le
\frac{1}{\sqrt N}
\sum_{i=1}^K
\sqrt{
G_i(d)\bigl(E-G_i(d)\bigr)
}
\\
&\le
\sqrt{\frac{K}{N}}
\left(
\sum_{i=1}^K
G_i(d)\bigl(E-G_i(d)\bigr)
\right)^{1/2}.
\end{align*}
Since
$
\sum_{i=1}^K G_i(d)=E,
$
we have
\[
\sum_{i=1}^K
G_i(d)\bigl(E-G_i(d)\bigr)
=
E^2-\sum_{i=1}^K G_i(d)^2
\le E^2,
\]
which proves the result.
\end{proof}

Let \(\tilde\rho\) be the reflector constructed from Algorithm \ref{AlgSP}, and define the discrete GOAS measure by
\begin{equation}
    \mu^{K,N,\epsilon}_{\mathrm{GOAS}}
    :=(Q\circ T_{\tilde\rho})_\sharp\mu_I=Q_\sharp\mu_{R_{\tilde\rho}},
    \label{eq:discrete_GOAS_measure}
\end{equation}
where $\mu_{R_{\tilde\rho}}$ is given by \eqref{eq:discrete_measure_bR}.

\begin{thm}[Reflector construction error]
\label{thm:discrete_GOAS_error}
Let $\hat d$ be the focal vector returned by Algorithm~\ref{AlgSP}, and let $\tilde\rho$ be the corresponding
discrete reflector. Assume that $Q:\Omega_O\to\Omega$ is $L_Q$-Lipschitz.
If~\eqref{eq:computed_stopping_criterion} holds, then
\begin{align}
\label{eq:discrete_GOAS_error}
W_p
\left(
\mu^{K,N,\epsilon}_{\mathrm{GOAS}},
\mu_t^K
\right)
\le
L_QD_O
\left(
\frac{
\eta_N(\hat d)+KE\epsilon
}{2}
\right)^{1/p}.
\end{align}
\end{thm}

\begin{proof}
For each $i$, the triangle inequality gives
\[
\left|
G_i(\hat d)-w_i
\right|
\le
\left|
G_i(\hat d)-\widetilde G_i(\hat d)
\right|
+
\left|
\widetilde G_i(\hat d)-w_i
\right|.
\]
Summing over $i$ yields
\begin{align}
\sum_{i=1}^K
\left|
G_i(\hat d)-w_i
\right|
&\le
\eta_N(\hat d)
+
\sum_{i=1}^K
\left|
\widetilde G_i(\hat d)-w_i
\right|.
\label{eq:mass_error_decomposition}
\end{align}
By the Cauchy--Schwarz inequality and
\eqref{eq:computed_stopping_criterion},
\begin{align*}
\sum_{i=1}^K
\left|
\widetilde G_i(\hat d)-w_i
\right|
&\le
\sqrt{
K
\sum_{i=1}^K
\left(
\widetilde G_i(\hat d)-w_i
\right)^2
}
\\
&\le
KE\epsilon.
\end{align*}
Hence,
\begin{align}
\sum_{i=1}^K
\left|
G_i(\hat d)-w_i
\right|
\le
\eta_N(\hat d)+KE\epsilon.
\label{eq:mass_total_error}
\end{align}

Since $\mu_{R_{\tilde\rho}}$ and $\mu_O^K$ are supported on the same points $\{y_i\}_{i=1}^K$, Lemma~\ref{lem:mass_perturbation} and \eqref{eq:mass_total_error} imply
\[
W_p
\left(
\mu_{R_{\tilde\rho}},
\mu_O^K
\right)
\le
D_O
\left(
\frac{
\eta_N(\hat d)+KE\epsilon
}{2}
\right)^{1/p}.
\]
Since
$
\mu_{\mathrm{GOAS}}^{K,N,\epsilon}
=
Q_\sharp\mu_{R_{\tilde\rho}}
$
and
$
\mu_t^K
=
Q_\sharp\mu_O^K,
$
from Lemma \ref{lem:Treg}, the Lipschitz push-forward estimate gives
\begin{align*}
W_p
\left(
\mu_{\mathrm{GOAS}}^{K,N,\epsilon},
\mu_t^K
\right)
&\le
L_Q
W_p
\left(
\mu_{R_{\tilde\rho}},
\mu_O^K
\right)
\\
&\le
L_QD_O
\left(
\frac{
\eta_N(\hat d)+KE\epsilon
}{2}
\right)^{1/p}.
\end{align*}
\end{proof}

\begin{cor}[Total Wasserstein error of discrete GOAS measure]\label{cor:dis_GOASmeasure}
Under the assumptions of Theorem~\ref{thm:discrete_GOAS_error},
\begin{align}
\mathbb E\,
W_p
\left(
\mu^{K,N,\epsilon}_{\mathrm{GOAS}},
\mu_t
\right)
\le
\delta_{K}
+
L_QD_O
\left[
\frac12
\left(
E\sqrt{\frac KN}
+
KE\epsilon
\right)
\right]^{1/p}.
\end{align}
\end{cor}

\begin{proof}
From the triangle inequality,
$$
W_p
\left(
\mu^{K,N,\epsilon}_{\mathrm{GOAS}},
\mu_t
\right)
\le
W_p
\left(
\mu^{K,N,\epsilon}_{\mathrm{GOAS}},
\mu_t^K
\right)
+W_p
\left(
\mu^K_t,
\mu_t
\right).
$$
The result follows from Theorem~\ref{thm:discrete_GOAS_error}, Lemma \ref{lem:MC_mass_error} and Jensen's inequality.
\end{proof}

\subsection{Softmin smoothing error}\label{sec:softmin_error}

We now compare the nonsmooth discrete reflector map with the smooth reflection map of Proposition \ref{prop:softmin_reflector_map}.
Let \(\tilde\rho\) be the reflector constructed from the supporting hyperparaboloid Algorithm \ref{AlgSP}, and let $\tilde\rho_\lambda$ be its softmin regularization.
Define the softmin-smoothed GOAS measure by
\begin{equation}
    \mu^{K,N,\epsilon,\lambda}_{\mathrm{GOAS}}
    :=(Q\circ T_{\tilde\rho_\lambda})_\sharp\mu_I.
    \label{eq:smooth_GOAS_measure}
\end{equation}


\begin{thm}[Wasserstein error of the softmin smoothing]\label{thm:softmin_Wp_error}
Assume that $\tilde\rho$, $\tilde\rho_\lambda$ and their first derivatives are uniformly bounded.
If $Q:\Omega_O\to\Omega$ is $L_Q$-Lipschitz, then there exists a constant $C>0$, independent of $\lambda$, such that
\begin{align}\label{eq:Wp_softmin_error_target}
W_p\bigl(\mu^{K,N,\epsilon,\lambda}_{\mathrm{GOAS}},\mu^{K,N,\epsilon}_{\mathrm{GOAS}}\bigr)
\le
CL_Q(\lambda e^{-\tau/\lambda}+e^{-\tau/\lambda}),
\end{align}
where $C$ is a positive constant.
\end{thm}

\begin{proof}
From \eqref{eq:softmin_reflection_map}, the smooth reflector map is given by
\[
T_{\tilde\rho_\lambda}(x)
=
x+
\frac{
2\tilde\rho_\lambda
\bigl(
D\tilde\rho_\lambda
-
x(\tilde\rho_\lambda+D\tilde\rho_\lambda\cdot x)
\bigr)
}{
\tilde\rho_\lambda^2
+
|D\tilde\rho_\lambda|^2
-
(D\tilde\rho_\lambda\cdot x)^2
}.
\]
The discrete reflector map $T_{\tilde\rho}$ is obtained by the same formula with $\tilde\rho_\lambda$ replaced by $\tilde\rho$.
Since involved quantities are uniformly bounded, the map
\[
(\rho,D\rho)\mapsto T_\rho
\]
is Lipschitz. Hence,
\[
|T_{\tilde\rho_\lambda}(x)-T_{\tilde\rho}(x)|
\le
C\left(
|\tilde\rho_\lambda(x)-\tilde\rho(x)|
+
|D\tilde\rho_\lambda(x)-D\tilde\rho(x)|
\right).
\]
By Theorem \ref{thm:rhoSmo}, under the unique-minimizer condition,
\[
|\tilde\rho_\lambda(x)-\tilde\rho(x)|
\le
C\lambda e^{-\tau/\lambda},
\qquad
|D\tilde\rho_\lambda(x)-D\tilde\rho(x)|
\le
C e^{-\tau/\lambda}.
\]
Therefore,
\begin{align}\label{eq:Terr}
|T_{\tilde\rho_\lambda}(x)-T_{\tilde\rho}(x)|
\le
C (\lambda e^{-\tau/\lambda}+e^{-\tau/\lambda}).
\end{align}

Define the coupling
\[
\gamma=(T_{\tilde\rho_\lambda},T_{\tilde\rho})_\sharp\mu_I .
\]
Then $\gamma$ has marginals $(T_{\tilde\rho_\lambda})_\sharp\mu_I$ and $(T_{\tilde\rho})_\sharp\mu_I$. Hence, by the definition of $W_p$,
\[
W_p^p\bigl((T_{\tilde\rho_\lambda})_\sharp\mu_I,(T_{\tilde\rho})_\sharp\mu_I\bigr)
\le
\int_{\Omega_I}
d_{S^n}\bigl(T_{\tilde\rho_\lambda}(x),T_{\tilde\rho}(x)\bigr)^p
\diff\mu_I(x).
\]
Moreover, for any $a,b\in S^n$,
\[
|a-b|^2=2(1-a\cdot b)
=2(1-\cos d_{S^n}(a,b))
=4\sin^2\left(\frac{d_{S^n}(a,b)}{2}\right).
\]
Since $0\le d_{S^n}(a,b)\le\pi$, we have
\[
d_{S^n}(a,b)\le \frac{\pi}{2}|a-b|.
\]
Therefore,
\[
W_p^p\bigl((T_{\tilde\rho_\lambda})_\sharp\mu_I,(T_{\tilde\rho})_\sharp\mu_I\bigr)
\le
\left(\frac{\pi}{2}\right)^p
\int_{\Omega_I}
|T_{\tilde\rho_\lambda}(x)-T_{\tilde\rho}(x)|^p
\diff\mu_I(x).
\]
Using \eqref{eq:Terr}, we obtain
\[
W_p\bigl((T_{\tilde\rho_\lambda})_\sharp\mu_I,(T_{\tilde\rho})_\sharp\mu_I\bigr)
\le
C (\lambda e^{-\tau/\lambda}+e^{-\tau/\lambda}).
\]

Finally, since $Q$ is Lipschitz, from Lemma \ref{lem:Treg}
\[
W_p\bigl((Q\circ T_{\tilde\rho_\lambda})_\sharp\mu_I,(Q\circ T_{\tilde\rho})_\sharp\mu_I\bigr)
\le
L_Q
W_p\bigl((T_{\tilde\rho_\lambda})_\sharp\mu_I,(T_{\tilde\rho})_\sharp\mu_I\bigr),
\]
which proves \eqref{eq:Wp_softmin_error_target}.
\end{proof}

\begin{cor}[Total Wasserstein error of smooth GOAS measure]\label{cor:smo_GOASmeasure}
Under the assumptions of Theorem~\ref{thm:discrete_GOAS_error} and Theorem~\ref{thm:softmin_Wp_error},
\begin{align*}
\mathbb E\,
W_p
\left(
\mu^{K,N,\epsilon,\lambda}_{\mathrm{GOAS}},
\mu_t
\right)
&\le
\delta_{K}
+
L_QD_O
\left[
\frac12
\left(
E\sqrt{\frac KN}
+
KE\epsilon
\right)
\right]^{1/p}\\
&+CL_Q(\lambda e^{-\tau/\lambda}+e^{-\tau/\lambda}).
\end{align*}
\end{cor}

\begin{proof}
From the triangle inequality,
$$
W_p
\left(
\mu^{K,N,\epsilon,\lambda}_{\mathrm{GOAS}},
\mu_t
\right)
\le
W_p
\left(
\mu^{K,N,\epsilon,\lambda}_{\mathrm{GOAS}},
\mu^{K,N,\epsilon}_{\mathrm{GOAS}}
\right)
+W_p
\left(
\mu^{K,N,\epsilon}_{\mathrm{GOAS}},
\mu_t^K
\right)
+W_p
\left(
\mu^K_t,
\mu_t
\right).
$$
The result follows from Theorem~\ref{thm:softmin_Wp_error}, Theorem~\ref{thm:discrete_GOAS_error}, Lemma \ref{lem:MC_mass_error} and Jensen's inequality.
\end{proof}

\begin{rem}
The estimate in Corollary \ref{cor:smo_GOASmeasure} separates the principal approximation mechanisms of the numerical GOAS. The term $\delta_{K}$ represents the target-measure approximation and therefore accommodates both density-based and sample-based target discretizations. The quantity $\eta_{N}$ represents the Monte Carlo evaluation of the reflector-induced masses, while $\epsilon$ measures the finite accuracy of the supporting hyperparaboloid construction. The final term is the additional error introduced by softmin smoothing.

No additional error associated with the local proposal used in the primal--dual consistency resampling procedure is considered here.
\end{rem}

\section{Numerical experiments}\label{sec:numexp}

In this section, we investigate the numerical performance of the proposed far-field geometric optics approximation sampling (GOAS) framework through a series of representative examples. The experiments are designed to examine
(i) sampling from strongly non-Gaussian distributions,
(ii) comparison with Markov chain Monte Carlo (MCMC) and transport-map-based samplers,
(iii) stability with respect to observational noise and target discretization,
(iv) Wasserstein error behavior predicted by the analysis in Section~\ref{sec:wasserstein_error}, and
(v) applicability to nonlinear PDE-constrained Bayesian inverse problems provided in Appendix \ref{Appendix:DetailsNumExpADR}.
The Bayesian formulation is given in Appendix~\ref{Appendix:BayInvPro}.

In all experiments, samples are generated using the primal--dual consistency resampling strategy introduced in Algorithm~\ref{alg:dualresim}. This choice allows us to work directly with the nonsmooth piecewise reflector produced by the supporting hyperparaboloid construction, while preserving the primal--dual cell correspondence without introducing an additional smoothing parameter. The alternative softmin-based realization provides an explicit smooth transport map, but is not repeated numerically here, since an analogous smooth-reflector sampling strategy has already been investigated in our near-field GOAS study. The present experiments therefore focus on the primal--dual realization that is specific to the far-field formulation.

The reflector surfaces are constructed using Algorithm~\ref{AlgSP} with tolerance parameter $\epsilon=10^{-4}$, and we use the density-based target discretization.
Let
$
y=(y_1,\ldots,y_{n+1})\in S^n,
z=(z_1,\ldots,z_n)\in\mathbb R^n.
$
Throughout the numerical experiments, we use the stereographic projection
\[
Q(y)
=
\frac{(y_1,\ldots,y_n)}{1-y_{n+1}},
\]
whose inverse is
\[
Q^{-1}(z)
=
\frac{(2z_1,\ldots,2z_n,|z|^2-1)}
{|z|^2+1}.
\]
Thus,
\[
Q:S^n\setminus\{(0,\ldots,0,1)\}
\longrightarrow\mathbb R^n
\]
is a diffeomorphism.

For all experiments, the source distribution is taken to be uniform on the north hemisphere.
Source rays are generated by normalizing Gaussian random vectors,
\[
x=\frac{X}{|X|},
\qquad
X\sim\mathcal N(0,I),
\]
and set $x_{n+1}=|x_{n+1}|$ where $I$ denotes the identity matrix in $\mathbb R^{n+1}$.

\subsection{GOAS for strongly non-Gaussian distributions}
\label{subsec:GOASvsMCMCs}

In this section, we first compare GOAS with several classical MCMC methods, including the Metropolis--Hastings (MH) algorithm \cite{robert2004}, slice sampling \cite{neal2003}, Hamiltonian Monte Carlo (HMC) \cite{duane1987,neal2011}, and the Metropolis-adjusted Langevin algorithm (MALA) \cite{roberts1998}.
The comparison is performed on several representative two-dimensional non-Gaussian target distributions, including Funnel, Banana, Mixture of Gaussians (MoG), Ring, and Cosine distributions \cite{wenliang2019,jaini2019}.
These examples exhibit a range of geometric complexities and multimodality.

\begin{figure}[htbp]
 \centering
 \includegraphics[scale=0.45]{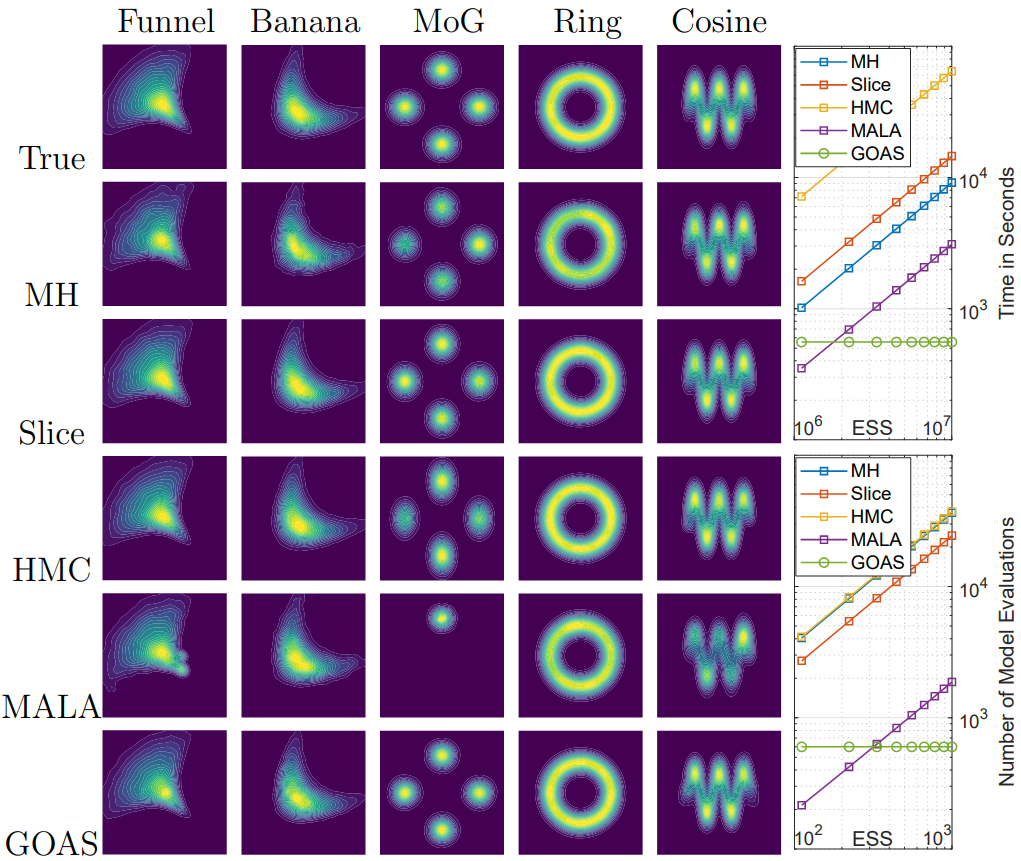}
 \caption{
Comparison of GOAS and MCMCs for sampling non-Gaussian target distributions.
The first five columns display the true densities together with the kernel density estimates obtained from GOAS and MCMCs.
For the MoG example, the rightmost column shows the computational time and the number of target density evaluations versus the effective sample size (ESS).
}
\label{fig:GOASvsMCMCs}
\end{figure}

In GOAS, the target measures are discretized using Hammersley low-discrepancy sequences together with the density-based discretization.
The true densities and the kernel density estimates obtained by GOAS and the MCMC methods are shown in Figure \ref{fig:GOASvsMCMCs}.
For the MoG example, the rightmost column additionally reports the computational time and the number of target evaluations with respect to the effective sample size (ESS); see Appendix \ref{AppendixESSHD} for details of the ESS computation.

The numerical results show that GOAS captures the geometric structures of all target distributions while directly generating independent samples through reflector-induced transport.
Among the MCMC methods, slice sampling also provides accurate density approximations for most examples, including the multimodal MoG distribution.
However, for MCMC methods, both the computational cost and the number of target density evaluations increase with the required ESS due to sample correlations.
In contrast, GOAS generates independent samples once the reflector surface has been constructed, and its sampling stage does not suffer from correlation-induced ESS limitations.
This behavior is particularly attractive when a large number of samples is required, especially in Bayesian inverse problems involving computationally expensive forward models.

\subsection{Comparison with polynomial transport maps}\label{subsec:GOASvsTM}

We next compare GOAS with a polynomial transport-map (TM) sampler on two Bayesian inverse problems: the Biochemical Oxygen Demand (BOD) model \cite{sullivan2010,marzouk2016} and the Euler--Bernoulli beam problem \cite{peherstorfer2019}. The two examples are chosen so that their posterior distributions exhibit similar banana-shaped geometries, while their forward models have markedly different computational complexity. This setting allows us to compare the approximation quality and computational cost of the two sampling strategies under similar posterior geometry.

A key distinction between the implementations considered here is that GOAS requires only evaluations of the posterior density at the target discretization points, whereas the TM construction involves an iterative optimization procedure using posterior-density and gradient information. For the BOD problem, the forward model is explicit and the required derivatives can be evaluated analytically. For the Euler--Bernoulli beam problem, each posterior evaluation requires the solution of a differential equation, and the gradients used in the TM construction are approximated by finite differences.
The TM implementation follows \cite{baptista2023,parno2022}.


\subsubsection*{Biochemical Oxygen Demand Model}\label{subsubsec:BODModel}
For the Biochemical Oxygen Demand (BOD) problem, the objective is to estimate two coefficients in a time-dependent model of oxygen demand, which serves as an indicator of biological activity in a water sample.
The simple time-dependent forward model is given by
$$
\mathcal{G}(t)=A(1-\exp(-Bt)),
$$
where $A=[0.4+0.4(1+erf(\theta_1/\sqrt{2}))]$ and $B=[0.01+0.15(1+erf(\theta_2/\sqrt{2}))]$.
The measurement data is obtained by $y(t)=\mathcal{G}(t)+\eta$ where $\eta\sim \mathcal{N}(0,10^{-3})$.
The objective is to characterize the posterior density of parameters $\theta=(\theta_1,\theta_2)$ knowing observation of the system at time $t={1,2,3,4,5}$, i.e., $y=(y(1),y(2),y(3),y(4),y(5))=(0.18,0.32,0.42,0.49,0.54)$.
Using a standard Gaussian prior for parameters $\theta$ and from equation \eqref{PosDen}, the posterior density is given by
$$
\pi(\theta)=\exp\biggl[-\frac{1}{2}\frac{1}{10^{-3}}\sum_{i=1}^5(\mathcal{G}(t_i)-y(t_i))^2-\frac{1}{2}\theta\theta^{T}\biggr].
$$
Obviously, it is easy to obtain the gradients of the posterior density.

\subsubsection*{Euler Bernoulli beam problem}\label{subsubsec:EBbProblem}
Consider a cantilever beam of $L>0$ length modeled by the $\Omega=[0,L]$, and the left boundary is fixed at the origin and the right boundary is free.
The Euler-Bernoulli beam is governed by the fourth-order differential equation
\begin{align}\label{EBbProblemPDE}
\frac{\partial^2}{\partial x^2}\biggl(E(x)\frac{\partial^2}{\partial x^2}u(x)\biggr)=f(x),\quad x\in\Omega,
\end{align}
where $u:\Omega\to\mathbb{R}$ is the vertical deflection of the beam and $f:\Omega\to\mathbb{R}$ is the load.
The effective stiffness of the beam is given by $E:\Omega\to\mathbb{R}$ describes beam geometry and material properties.
We consider a smoothed piecewise constant approximation $\widetilde{E}_2:\Omega\times \mathbb{R}\to \mathbb{R}$ of the stiffness $E$ that depends on parameters $\theta=(\theta_1,\theta_2)$.
The forward operator $\mathcal{G}$  is then given by numerically solving \eqref{EBbProblemPDE} with stiffness $\widetilde{E}_2$.
The Euler Bernoulli beam problem is to infer the parameters $\theta$ in $\widetilde{E}_2$ from the measurements of displacement $u$ in $\Omega$.
In order to compare the effects of the different problems on the GOAS and TM methods, we set the appropriate parameters to achieve a consistent geometric structure of the posterior distribution, specifically the ``banana shape".
For details on the inverse problem setup, refer to Appendix \ref{Appendix:DetailsNumExpEBb}.
\begin{figure}[htbp]
  \centering
  \subfloat[]
  {   \includegraphics[width=0.2\textwidth,height=0.2\textwidth]{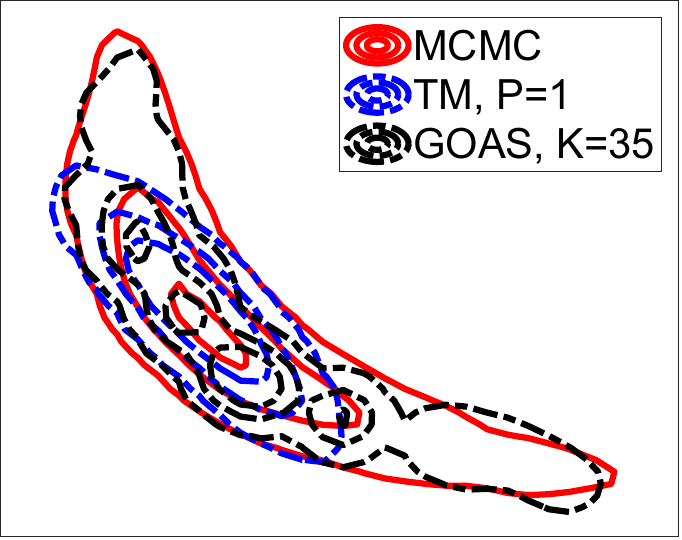}
  }
  \subfloat[]
  {   \includegraphics[width=0.2\textwidth,height=0.2\textwidth]{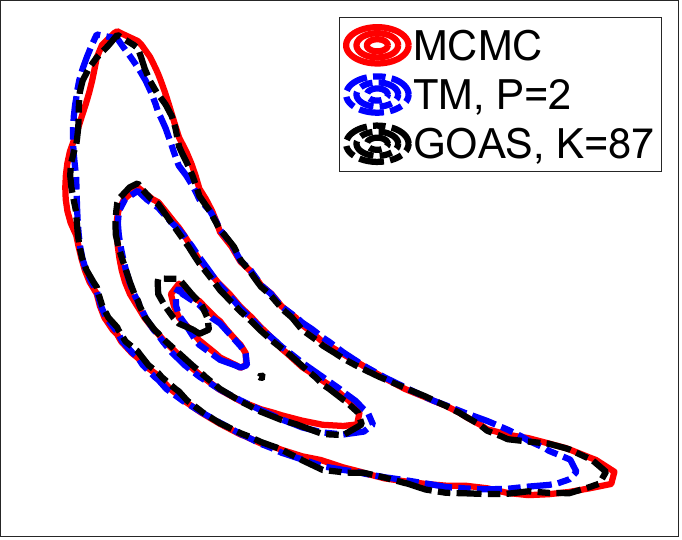}
  }
  \subfloat[]
  {   \includegraphics[width=0.2\textwidth,height=0.2\textwidth]{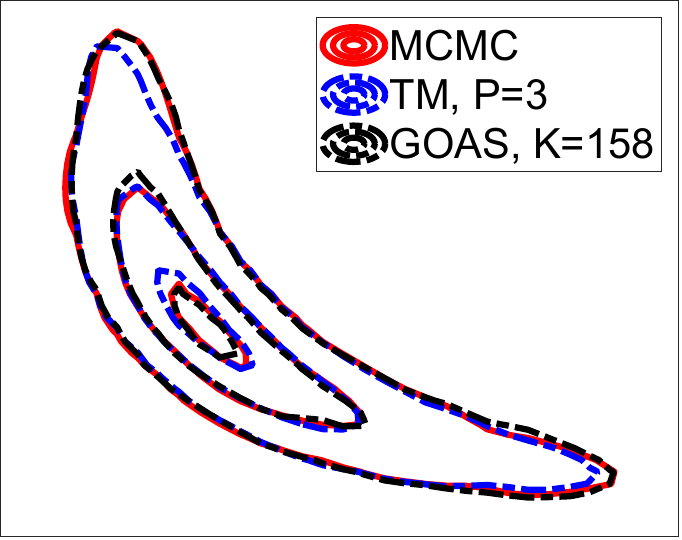}
  }\\
  \subfloat[]
  {   \includegraphics[width=0.2\textwidth,height=0.2\textwidth]{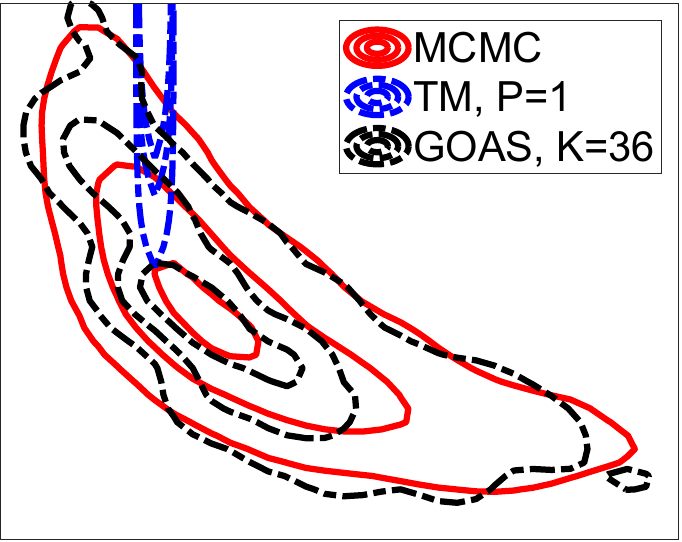}
  }
  \subfloat[]
  {   \includegraphics[width=0.2\textwidth,height=0.2\textwidth]{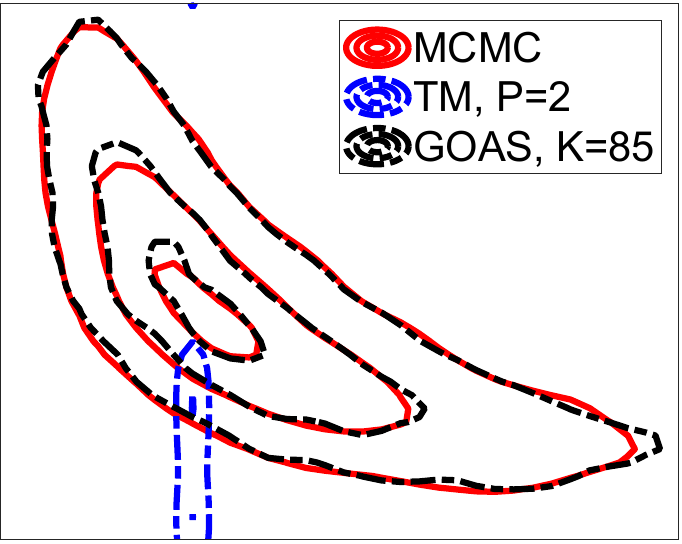}
  }
  \subfloat[]
  {   \includegraphics[width=0.2\textwidth,height=0.2\textwidth]{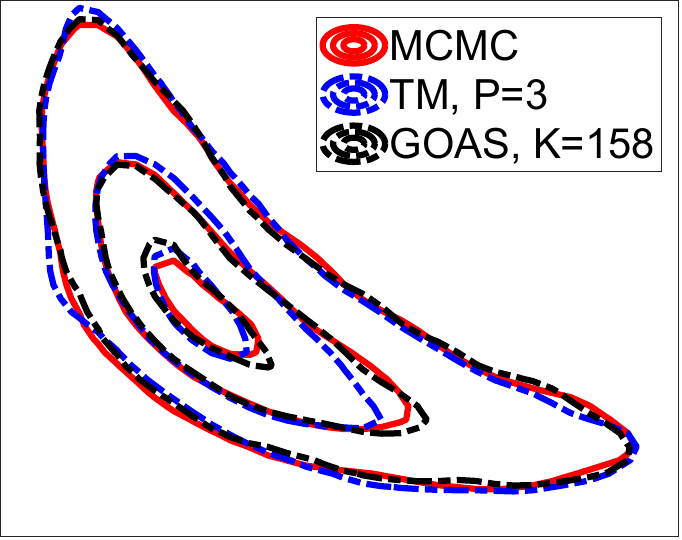}
  }
  \caption{
Kernel density estimates obtained by GOAS and the polynomial transport-map (TM) sampler for the BOD model (top row) and the Euler--Bernoulli beam problem (bottom row).
GOAS uses $K=35,87,158$ paraboloids, while TM uses total-order Hermite polynomial maps of degree $P=1,2,3$.
The MCMC estimate is shown as the reference posterior.
}
 \label{fig:GOASvsTM1}
\end{figure}



\begin{table}[htbp]
\centering
\caption{
Sample mean, variance, skewness, and kurtosis of the posterior parameters $\theta=(\theta_1,\theta_2)$ obtained by GOAS and TM methods for the BOD model and the Euler--Bernoulli beam problem.
}
\begin{tabular}{c|c|cccc}
\hline
\multicolumn{2}{c}{Method}                         & Mean & Variance & Skewness & Kurtosis \\ \hline
\multicolumn{6}{c}{\cellcolor{blue!10}Biochemical Oxygen Demand model} \\ \hline
\rowcolor{blue!10}
\multicolumn{2}{c|}{\cellcolor{blue!10}MCMC ``true"} & (0.0340, 0.9224) &(0.1443, 0.3716)          &(1.781, 0.6741)  &(7.329, 3.440)          \\ \hline
\multirow{3}{*}{TM}           & P=1   &(0.0921, 0.9682) &(0.0289, 0.1155) &(-0.0040, -0.0654)&(2.884, 2.926)          \\
                                                & P=2   &(0.0328, 0.9111) &(0.1264, 0.4004) &(1.609, 1.156) &(7.373, 5.784 )          \\
                                                & P=3   &(0.0383, 0.9147)      &(0.1611, 0.3831)          &(2.737, 0.9659)          &(20.88, 5.115)          \\ \hline
\rowcolor{gray!10}
\multicolumn{1}{c|}{\cellcolor{gray!10}} & K=35   &(0.0609, 0.8384)      &(0.1697, 0.3715)          &(1.858, 0.8179)          &(8.213, 3.611)          \\
\rowcolor{gray!10}
\multicolumn{1}{c|}{\cellcolor{gray!10}}& K=87   &(0.0467, 0.9014)      &(0.1509, 0.3854)          &(1.810, 0.6721)         &(7.901, 3.296)          \\
\rowcolor{gray!10}
\multicolumn{1}{c|}{\multirow{-3}{*}{\cellcolor{gray!10}GOAS}} & K=158   &(0.0459, 0.9302)      &(0.1704, 0.4019)          &(1.919, 0.5989)          &(7.853, 3.184)            \\ \hline

\multicolumn{6}{c}{\cellcolor{blue!10}Euler Bernoulli Beam Problem} \\ \hline
\rowcolor{blue!10}
\multicolumn{2}{c|}{\cellcolor{blue!10}MCMC ``true"} &(0.7316, 1.397)      &(0.0338, 0.1025)          &(1.126, 0.6622)          &(4.498, 3.310)          \\ \hline
\multirow{3}{*}{TM}           & P=1   &(0.5793, 3.410)      & (0.0001, 1.040)         &(0.0079, -0.0060)          &  (3.053, 3.057)        \\
                                                & P=2   &(0.6584, 1.755)  &(0.0000, 0.2013)   & (-11.31, 0.5504)         &(217.2, 5.807)          \\
                                                & P=3   &(0.7204, 1.423)      &(0.0365, 0.1090)          &(1.797, 0.5487)          & (9.847, 3.146)         \\ \hline
\rowcolor{gray!10}
\multicolumn{1}{c|}{\cellcolor{gray!10}} & K=36   &(0.7218, 1.403)      &(0.0321, 0.1144)          &(1.219, 0.5752)          &(5.135, 2.902)         \\
\rowcolor{gray!10}
\multicolumn{1}{c|}{\cellcolor{gray!10}}& K=85   &(0.7314, 1.407)      &(0.0346, 0.1062)          &(1.165, 0.6730)          &(4.673, 3.232)          \\  
\rowcolor{gray!10}
\multicolumn{1}{c|}{\multirow{-3}{*}{\cellcolor{gray!10}GOAS}} & K=158   &(0.7303, 1.405)      &(0.0344, 0.1066)          &(1.159, 0.7187)          &(4.672, 3.437)            \\ \hline
\end{tabular}
\label{tab:GOASvsTM1}
\end{table}

\begin{table}[htbp]
\centering
\caption{
Offline computational cost, online sampling time, and number of forward model evaluations required by GOAS and TM methods for the BOD model and the Euler--Bernoulli beam problem at the same effective sample size.
}
\begin{tabular}{ccccccccc}
\hline
\multirow{2}{*}{Method} &
\multicolumn{2}{c}{Offline}& \multirow{2}{*}{Online time (s)}&
\multirow{2}{*}{Total (h)} \\
\cline{2-3} 
&model evaluation & time (s) \\ \hline
\multicolumn{5}{c}{\cellcolor{blue!10}Biochemical Oxygen Demand model} \\ \hline
TM, P=3 & $7.7\times10^5$ & 6.902 & 0.0213 & 0.0019  \\ \hline
\rowcolor{gray!10}
GOAS, K=158 & $1.16\times10^4$ & 165.4 & 0.0653 & 0.0459  \\ \hline
\multicolumn{5}{c}{\cellcolor{blue!10}Euler Bernoulli Beam Problem} \\ \hline
TM, P=3 & $9.48\times 10^6$ & $9.8078\times 10^4$ &  0.4096 & 27.24  \\ \hline
\rowcolor{gray!10}
GOAS, K=158 & $1.1\times10^4$ & 157.5  & 0.0808 & 0.0437  \\ \hline
\end{tabular}
\label{tab:GOASvsTM2}
\end{table}

\subsubsection*{Comparison results}

For GOAS, a bounded computational target region is first identified using a pilot set of prior samples together with a threshold on the posterior density. The target measure is then approximated over this region using Hammersley low-discrepancy points, with weights determined by the density-based discretization in \eqref{eq:density_based_weights}. In both examples, $10^4$ prior samples are used to identify the computational region. The threshold is chosen so that the retained region contains the dominant posterior mass.
For TM methods, total-order Hermite polynomial maps of degree $P=1,2,3$ are employed.
For the Euler--Bernoulli beam problem, second-order central finite differences are used to approximate the gradients of the posterior density.

Figure \ref{fig:GOASvsTM1} shows the kernel density estimates obtained by GOAS with $K=35,87,158$ hyperparaboloids and TM methods with polynomial degrees $P=1,2,3$.
The MCMC results are used as reference densities.

As the number of hyperparaboloids increases, the GOAS approximations progressively approach the reference posterior in both examples. In particular, comparable values of $K$ provide similar reconstruction quality for the two banana-shaped posterior distributions. This observation is consistent with the geometric nature of the GOAS construction: once the target density has been evaluated at the discretization points, the reflector approximation is determined primarily by the resulting target-measure geometry rather than by the computational complexity of the forward model.

For the polynomial TM sampler, increasing the polynomial degree generally improves the approximation. In the BOD example, the quadratic and cubic maps reproduce the reference posterior reasonably well. In the Euler--Bernoulli example, the low-order maps show larger discrepancies, whereas the cubic map provides a substantially improved approximation. These differences may reflect both the limited expressiveness of low-order polynomial maps and the additional numerical error introduced by the finite-difference gradient evaluations.

Table~\ref{tab:GOASvsTM1} reports the sample mean, variance, skewness, and kurtosis for the two posterior parameters.
As $K$ increases, the GOAS estimates generally approach the MCMC reference values for both inverse problems, in agreement with the density estimates in Figure~\ref{fig:GOASvsTM1}.
The posterior means and variances are already captured reasonably well at moderate values of $K$, while the higher-order moments exhibit somewhat larger fluctuations.
For the TM sampler, increasing the polynomial degree generally improves the lower-order statistics. The higher-order moments are more sensitive to the map approximation, particularly for the Euler--Bernoulli example.

Finally, Table~\ref{tab:GOASvsTM2} compares the offline computational cost, online sampling time, and number of forward-model evaluations for GOAS with $K=158$ and the cubic TM sampler, using the same ESS.

Both methods have a relatively small online sampling cost once their offline constructions have been completed. Their main difference therefore lies in the offline stage. GOAS constructs the reflector from posterior density evaluations at prescribed target points, whereas the TM implementation considered here solves an optimization problem involving both posterior values and gradients.

For the BOD problem, where the forward model and its derivatives are inexpensive, the cubic TM construction is faster than GOAS. The situation changes substantially for the Euler--Bernoulli problem. Here, the finite-difference gradient evaluations required by the TM optimization lead to approximately $9.48\times10^6$ forward solves, compared with approximately $1.1\times10^4$ forward solves for GOAS. Accordingly, the total computational time is about $27.24$ hours for TM and $0.0437$ hours for GOAS in this experiment.

These results illustrate a potential advantage of GOAS for PDE-constrained Bayesian inverse problems: its reflector construction is gradient-free and requires only posterior-density evaluations at the target discretization points. This can substantially reduce the offline cost when each forward-model evaluation is expensive.

\begin{figure}[htbp]
  \centering
  \subfloat[Noise Level=10\%]
  {   \includegraphics[width=0.3\textwidth,height=0.25\textwidth]{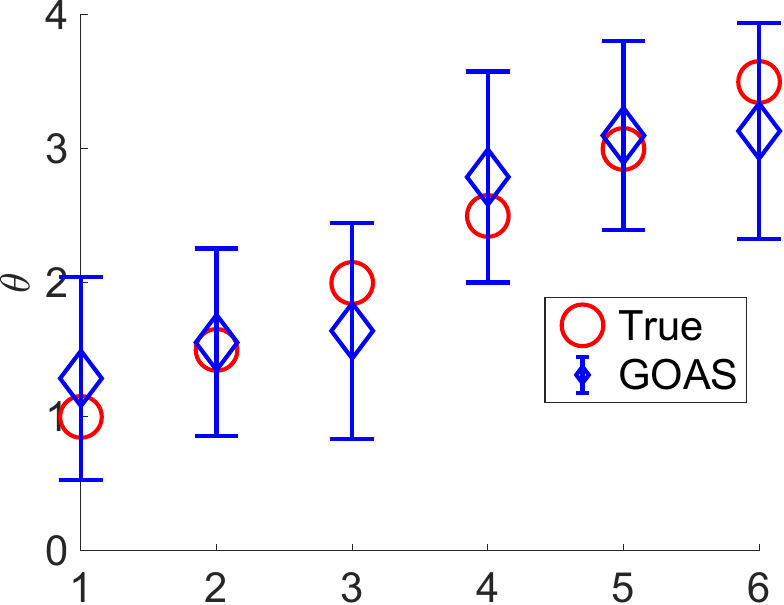}
  }
  \subfloat[Noise Level=5\%]
  {   \includegraphics[width=0.3\textwidth,height=0.25\textwidth]{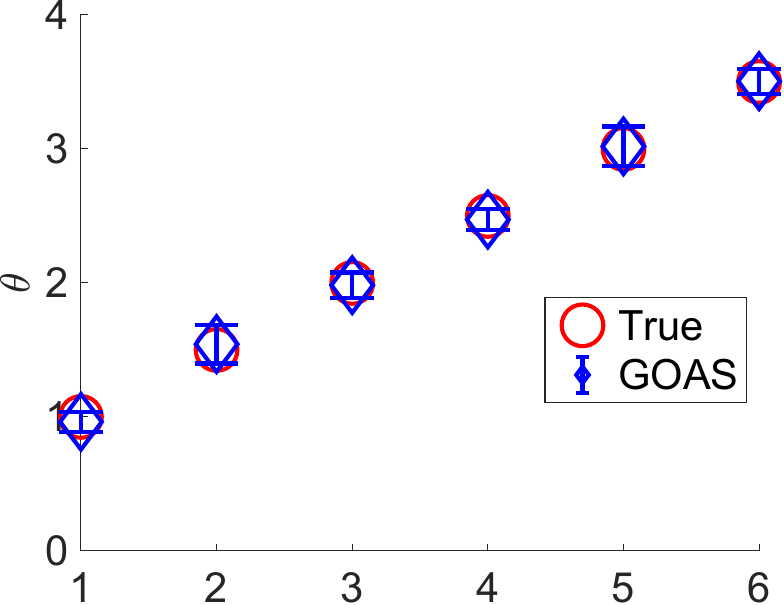}
  }
  \subfloat[Noise Level=1\%]
  {   \includegraphics[width=0.3\textwidth,height=0.25\textwidth]{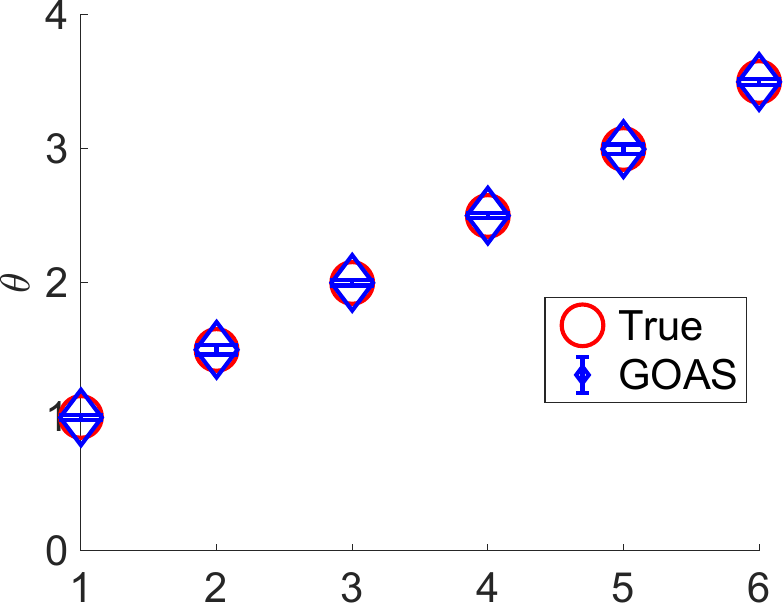}
  }
  \caption{
Posterior means and corresponding $95\%$ credible intervals of the parameters $\theta$ obtained by GOAS with $K= 800$ for the acoustic source localization under noise levels of $10\%$ (a), $5\%$ (b), and $1\%$ (c).
}
 \label{fig:exa11}
\end{figure}



\begin{figure}[htbp]
  \centering
   \subfloat[low-discrepancy-like design]
  {   \includegraphics[width=0.45\textwidth,height=0.3\textwidth]{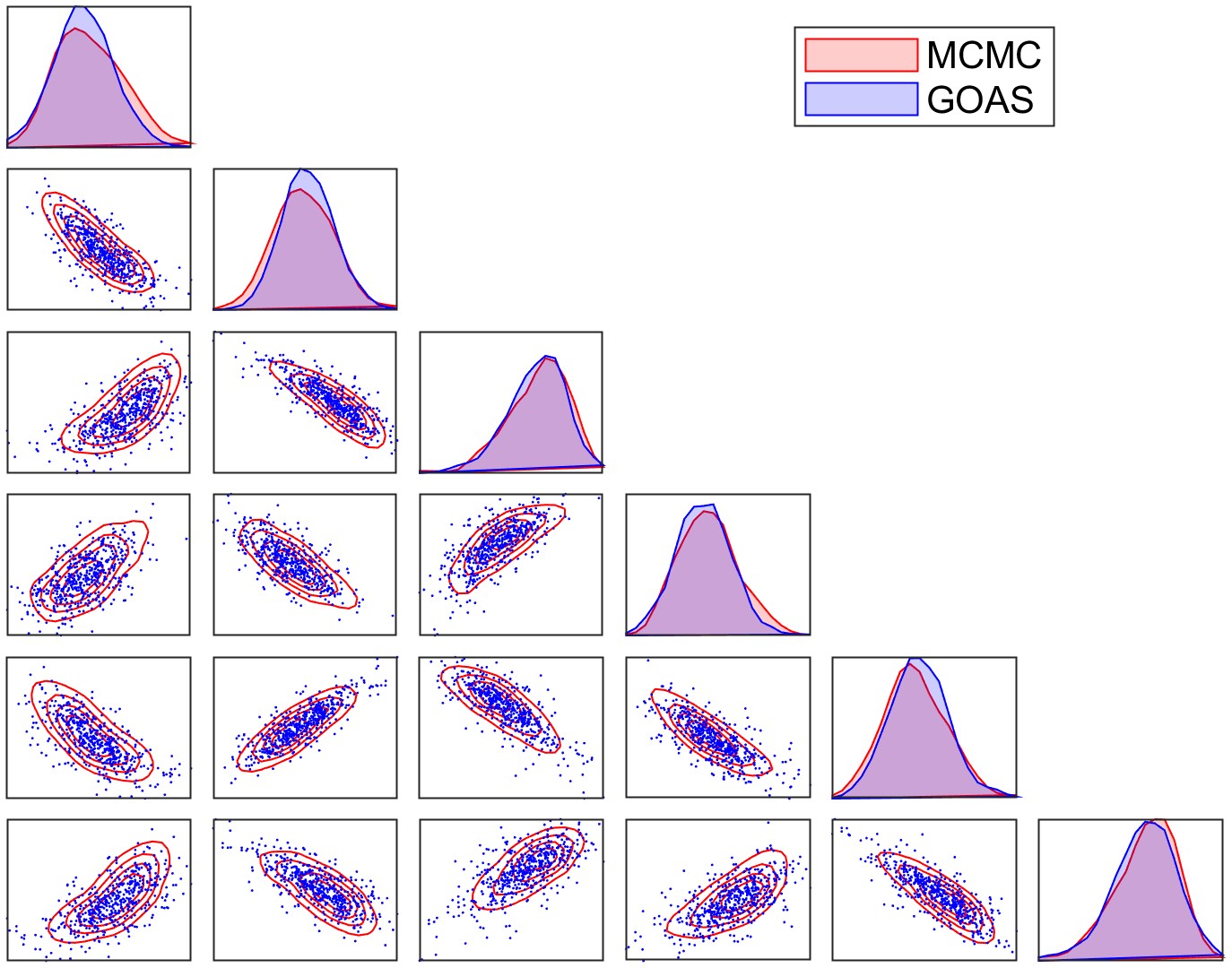}
  }
  \subfloat[importance-sampling strategy]
  {   \includegraphics[width=0.45\textwidth,height=0.3\textwidth]{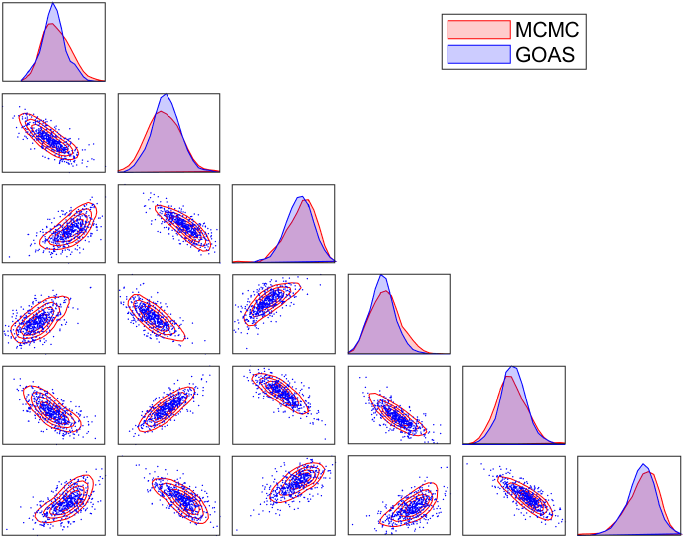}
  }
\caption{
Posterior samples and marginal density estimates obtained by GOAS using a low-discrepancy-like discretization (a) and an importance-sampling strategy discretization (b) of the target measure for the acoustic source localization at a $5\%$ noise level.
The MCMC results are shown as reference distributions.
}
\label{fig:exa12}
\end{figure}

\subsection{Stability validation: acoustic source localization}\label{subsec:exa1}

We next investigate the stability of GOAS with respect to the discretization of the target measure through an inverse acoustic source problem. Consider the Helmholtz equation
\begin{align}
\label{Hel}
\Delta u^s+k^2u^s=S
\qquad\text{in }\mathbb R^2,
\end{align}
where $k>0$ is the wave number, $u^s$ is the scattered field, and
\begin{align}
\label{Exa1:source}
S=\sum_{i=1}^N\lambda_i\delta_{z_i}
\end{align}
is a collection of point sources with locations $z_i\in\Omega$ and
nonzero intensities $\lambda_i$.
The inverse problem is to recover the source locations from measurements
of the far-field pattern; see, e.g.,
\cite{liu2021,eller2009,el2011}.
Details of the numerical setup are given in Appendix~\ref{Appendix:DetailsNumExpLAS}.

We first examine the effect of observational noise. For each noise level, the posterior target measure is approximated by the density-based discretization described in \eqref{eq:density_based_weights} with low-discrepancy-like design. Figure~\ref{fig:exa11} shows the posterior means and corresponding $95\%$ credible intervals of the parameters $\theta=\{\theta_i\}_{i=1}^{6}$ obtained by GOAS under noise levels of $10\%$, $5\%$, and $1\%$. As the noise level decreases, the posterior distributions become more concentrated around the true parameter values. Correspondingly, the posterior means approach the ground truth and the credible intervals become narrower, reflecting the reduction in posterior uncertainty as the data become more informative.

We then examine the sensitivity of GOAS to the choice of target discretization.  We consider two density-based discretizations, i.e., low-discrepancy-like design and  importance-sampling strategy. Figure~\ref{fig:exa12} compares the posterior samples and marginal density estimates obtained from the two discretizations at a $5\%$ noise level. The MCMC results are used as reference distributions. Both point-set designs recover the dominant posterior structure and produce marginal distributions that agree well with the MCMC reference.

These results indicate that GOAS is stable with respect to different the density-based target discretization.
This observation is consistent with the GOAS measure stability analysis developed in Corollary~\ref{cor:GOA_measure_stability}.


\begin{figure}[htbp]
  \centering
  \subfloat[Noise Level=10\%]
  {   \includegraphics[width=0.22\textwidth,height=0.22\textwidth]{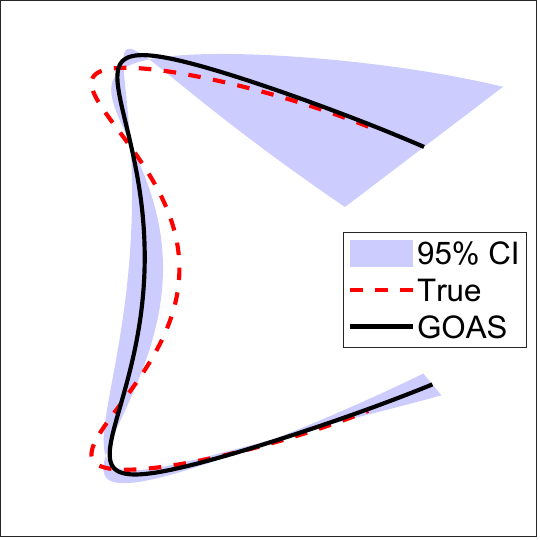}
  }
  \subfloat[Noise Level=5\%]
  {   \includegraphics[width=0.22\textwidth,height=0.22\textwidth]{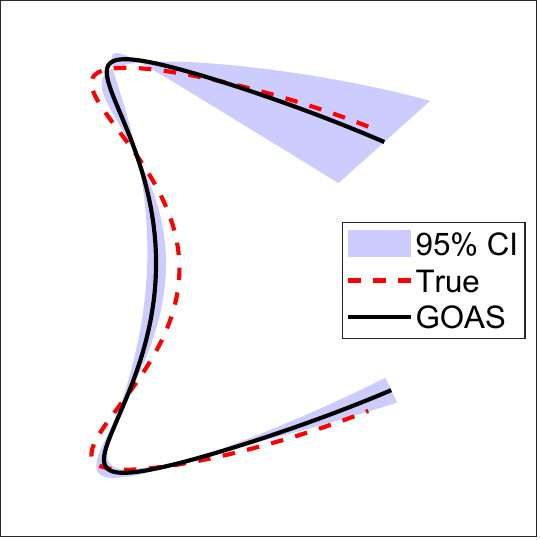}
  }
  \subfloat[Noise Level=1\%]
  {   \includegraphics[width=0.22\textwidth,height=0.22\textwidth]{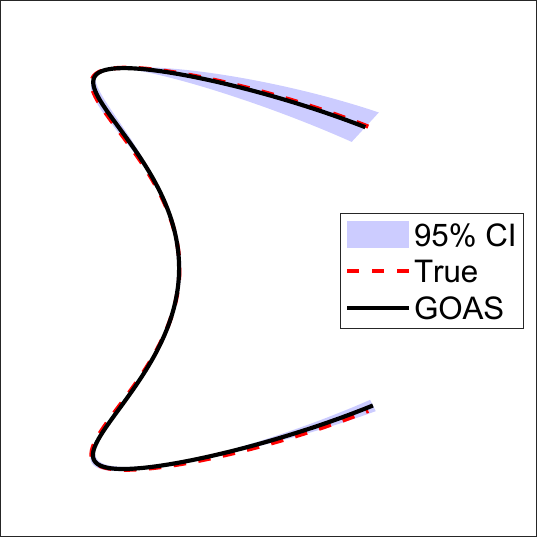}
  }
  \caption{
Posterior means and $95\%$ credible intervals of the parameter vector $\theta$ obtained by GOAS with approximately $K=890$ for the inverse scattering problem from an open arc at noise levels $10\%$ (a), $5\%$ (b), and $1\%$ (c).
}
 \label{fig:exa21}
\end{figure}

\begin{figure}[htbp]
  \centering
  \includegraphics[width=0.5\textwidth,height=0.3\textwidth]{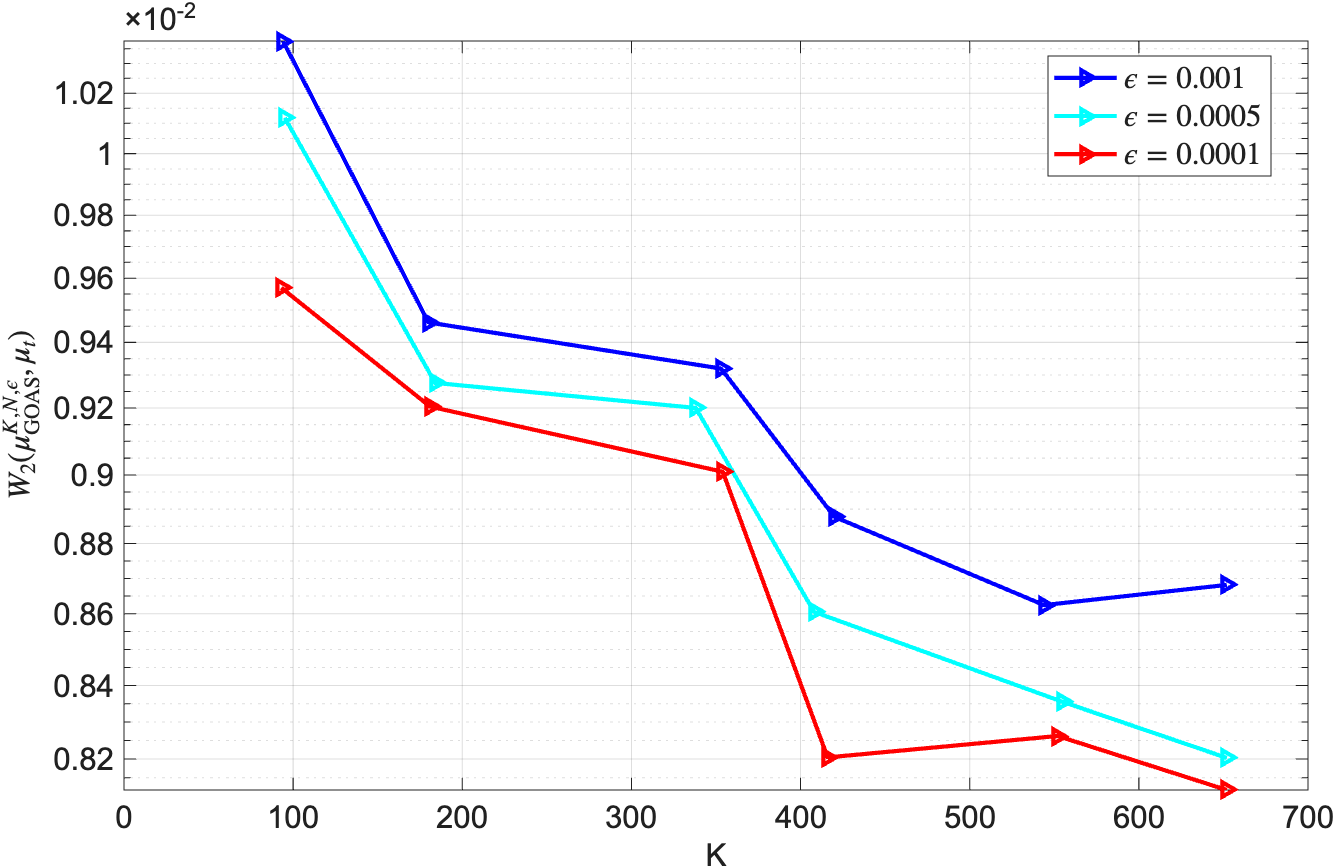}
\caption{Numerical Wasserstein discrepancy between reference posterior measure and discrete GOAS measure $W_2(\mu^{K,N,\epsilon}_{\mathrm{GOAS}},\mu_t)$ for reflector construction tolerances $\epsilon=0.001,0.0005,0.0001$ in \eqref{SPconErr} with respect to $K$ for inverse scattering from an open arc.}
\label{fig:exa22}
\end{figure}
\subsection{Wasserstein error validation: inverse scattering from an open arc}\label{subsec:exa2}
In this section, we examine the Wasserstein error behavior predicted by the analysis in Section~\ref{sec:wasserstein_error} through an inverse scattering problem for an open arc.
Consider the scattering of a time-harmonic acoustic wave by a thin sound-soft obstacle $\Gamma\subset\mathbb R^2$. The total field $u=u^i+u^s$ satisfies
\begin{align}
\label{IPopenarcPDE}
(\Delta+k^2)u &=0,\quad \text{in }\mathbb R^2\setminus\Gamma,\\
u&=0,\quad\text{on }\Gamma,
\end{align}
where $k>0$ is the wave number, $u^i$ is the prescribed incident field, and the scattered field $u^s$ satisfies the Sommerfeld radiation condition~\eqref{HelSRC}. The inverse problem is to recover the shape of the analytic arc $\Gamma$ from measurements of the far-field pattern corresponding to a single incident plane wave \cite{kress1995,kirsch2000}.
Further details of numerical setup are provided in Appendix~\ref{Appendix:DetailsNumExpIPopenarc}.

Figure~\ref{fig:exa21} first illustrates the posterior uncertainty quantified by GOAS. It shows the posterior means and corresponding $95\%$ credible intervals of the unknown parameter vector $\theta$ under noise levels of $10\%$, $5\%$, and $1\%$. As the noise level decreases, the posterior mean approaches the true parameter profile and the credible intervals become progressively narrower, reflecting the increasing information content of the data.

To investigate the convergence behavior predicted by the total Wasserstein error of discrete GOAS measure in Corollary \ref{cor:dis_GOASmeasure}, Figure \ref{fig:exa22} reports the Wasserstein error $W_2(\mu^{K,N,\epsilon}_{\mathrm{GOAS}},\mu_t)$ for reflector construction tolerances $\epsilon=0.001,0.0005,0.0001$ in \eqref{SPconErr}.
The MCMC results are considered as reference posterior $\mu_t$.
The Wasserstein distance is evaluated numerically using an entropically regularized optimal transport approximation with regularization parameter $10^{-3}$, computed by the Sinkhorn algorithm \cite{peyre2019computational}.

The results show that the Wasserstein discrepancy generally decreases as $K$ increases.
For a fixed $K$, decreasing the reflector tolerance $\epsilon$ also leads, in general, to a smaller discrepancy. These trends are consistent with the error decomposition in Corollary~\ref{cor:dis_GOASmeasure}: increasing the resolution of the target discretization reduces the target-approximation contribution, while a smaller stopping tolerance reduces the reflector-construction error.

\section{Conclusion}
\label{sec:conclusion}

In this work, we developed a far-field geometric optics approximation sampling (GOAS) framework for constructing transport-based samplers from prescribed target measures. The far-field reflector problem is naturally connected with optimal transport under the logarithmic cost, which provides the underlying transport structure of the method. In particular, reflector duality yields a natural primal--dual correspondence, so that, in addition to the forward transport from the source to the target, the associated dual reflector provides a natural mechanism for inverse transport.

For the numerical realization considered in this work, we use the supporting hyperparaboloid construction in general dimensions. Since this construction operates on a discrete target measure, both density-based and sample-based target representations arise naturally within the same framework. Moreover, the reflector is determined through target masses rather than derivatives of the target density, resulting in a gradient-free construction. The resulting discrete reflector is piecewise smooth, which leads to two complementary sampling realizations. The primal--dual consistency strategy works directly with the nonsmooth reflector by exploiting the reflector duality, whereas softmin regularization produces a smooth approximation from which an explicit transport map follows directly from the physical reflection law.

The supporting hyperparaboloid method is one numerical realization of the far-field GOAS framework rather than an essential restriction of the approach. Other numerical solvers for the corresponding far-field reflector or optimal transport problem may be incorporated while retaining the same underlying sampling principle.

We further established Wasserstein error estimates that separate the effects of target-measure discretization, Monte Carlo evaluation of reflector-induced masses, reflector construction, and smoothing. The numerical experiments demonstrate the ability of GOAS to represent strongly non-Gaussian distributions, its competitive performance relative to MCMC and polynomial transport-map methods, and its applicability to nonlinear PDE-constrained Bayesian inverse problems.

The particular supporting-hyperparaboloid realization adopted in this work has several numerical and analytical limitations. Since the construction is based on a discrete approximation of the target measure, the accuracy of the resulting sampler depends on the quality of this discretization. Moreover, the primal--dual consistency resampling procedure introduces a local within-cell approximation that is not included in the present Wasserstein analysis. These limitations are associated with the numerical realization used here rather than with the underlying far-field GOAS transport framework.
These issues, together with more efficient reflector construction, provide natural directions for future investigation.

\addcontentsline{toc}{section}{Acknowledgment}
\begin{ack}
The work described in this paper was supported by the NSF of China (12271151) and NSF of Hunan (2020JJ4166) and Postgraduate Scientific Research Innovation Project of Hunan Province (CX20240364).
\end{ack}


\begin{appendices}
The appendices includes mathematical formulation of far-field reflector problem, details on numerical experiments, additional numerical experiments, and all proof details.

\section{Far-field reflector geometry}\label{sec:farfield_geometry}
We first recall several standard notions from the far-field reflector; see, e.g., \cite{caffarelli1999,wang2004,caffarelli2008}.
For $y\in\Omega_O$ and focal parameter $d>0$, let $P(y)=\{xp_{y,d}(x): x\in S^n\}$ denote the paraboloid of revolution with focus at the origin, axial direction $y$, and polar radius
\[
p_{y,d}(x)=\frac{d}{1-x\cdot y},\qquad x\in S^n.
\]

\begin{defn}[Supporting paraboloid]\label{def:SP}
Let $R=R_\rho$ be a reflecting surface given by \eqref{RGraph}.
A paraboloid $P(y)$, $y\in\Omega_O$, is called supporting to $R_\rho$ at $x_0\rho(x_0)$ if
\[
\rho(x_0)=p_{y,d}(x_0),
\qquad
\rho(x)\le p_{y,d}(x),\quad \forall x\in\Omega_I.
\]
\end{defn}

\begin{defn}[Convex reflecting surface]\label{def:ConRefSur}
A reflecting surface $R=R_\rho$ is called convex with respect to $\Omega_I$ if at every point of its graph there exists a supporting paraboloid.
Let $B(y)$ denote the convex body bounded by the paraboloid $P(y)$.
If
\begin{align}\label{Rpolyhedron}
R_K=\partial B,
\qquad
B=\bigcap_{i=1}^K B(y_i),
\end{align}
where $y_i\in\Omega_O$ and $K\in\mathbb N^+$, then $R_K$ is called a polyhedral reflecting surface with respect to $\Omega_O$.
\end{defn}

Throughout this paper, unless otherwise stated, all reflecting surfaces are assumed to be convex.

\begin{defn}[Subdifferential]\label{DefSubdiff}
Let $R=R_\rho$ be a convex reflecting surface.
The subdifferential of $\rho$ is the set-valued map $T_\rho:\Omega_I\to\Omega_O$ defined by
\[
T_\rho(x)
=
\left\{
y\in\Omega_O:
\exists d>0 \ \text{such that } p_{y,d}
\text{ is supporting to } \rho \text{ at } x
\right\}.
\]
For any $E\subset\Omega_I$, we write
\[
T_\rho(E)=\bigcup_{x\in E}T_\rho(x).
\]
\end{defn}

For nonsmooth convex reflectors, the reflector-induced transport is described by the set-valued subdifferential map $T_\rho$. In the smooth case, the reflector induces a single-valued reflection map, which we denote by $T$.

\begin{rem}
The subdifferential map $T_\rho$ is single-valued at every differentiability point of $\rho$, and hence is single-valued almost everywhere for convex reflectors.
In general, $T_\rho$ may be multivalued for nonsmooth reflectors, whereas the smooth reflector formulation introduced later induces a single-valued reflection map.

Similarly, one defines the inverse visibility map $V:\Omega_O\to\Omega_I$ by
\[
V(y)
=
\left\{
x\in\Omega_I:
\exists d>0 \ \text{such that } p_{y,d}
\text{ is supporting to } \rho \text{ at } x
\right\}.
\]
For $w\subset\Omega_O$, the set
\[
V(w)=\bigcup_{y\in w}V(y)
\]
is called the visibility set of $w$.
If $R_\rho$ is smooth and $T_\rho$ is one-to-one, then $V=T_\rho^{-1}$.
\end{rem}

\begin{prop}\label{rhoTconvexRef}
Let $R_\rho$ be a convex reflecting surface.
Then its polar radius is given by
\begin{align}\label{ConRefRho}
\rho(x)=\inf_{y\in\Omega_O}p_{y,d(y)}(x),
\qquad x\in\Omega_I.
\end{align}
Moreover, at points where the minimizer is unique, the reflected direction is
\begin{align}\label{RefDirbySP}
T_\rho(x)
=
\arg\min_{y\in\Omega_O}p_{y,d(y)}(x).
\end{align}
\end{prop}

\begin{proof}
By convexity, for each $x\in\Omega_I$ there exists a supporting paraboloid touching $R_\rho$ at $x\rho(x)$.
Hence, for some $y\in\Omega_O$ and $d(y)>0$,
\[
\rho(x)=p_{y,d(y)}(x),
\qquad
\rho(x')\le p_{y,d(y)}(x'),
\quad \forall x'\in\Omega_I.
\]
Taking the infimum over all supporting paraboloids gives \eqref{ConRefRho}.
By the law of reflection, the axial direction of the supporting paraboloid at $x\rho(x)$ is precisely the reflected direction, which gives \eqref{RefDirbySP}.
\end{proof}

\section{Bayesian inverse problem framework}\label{Appendix:BayInvPro}
Let $\mathbf{X}$ be a separable Hilbert space, equipped with the Borel $\sigma$-algebra, and let
 $\mathcal{G}:\mathbf{X}\to \mathbb{R}^n$ be a measurable function known as the forward operator, which represents the connection between parameter and data in the mathematical model.
We aim to solve the inverse problems of finding the unknown model parameters $u$ in set $\mathbf{X}$ from measurement data $y\in\mathbb{R}^n$, which is usually generated by
\begin{align}\label{ForMod}
y=\mathcal{G}(u)+\eta,
\end{align}
where the noise $\eta$ is assumed to be a $n$-dimensional zero-mean Gaussian random variable with covariance matrix $\Sigma_{\eta}$.
In Bayesian inverse problems, the unknown model input and the measurement data are typically regarded as random variables.
Let $\mu_0$ be a prior probability measure on $\mathbb{R}^n$ and let $\Phi_y: \mathbb{R}^n\to [0,+\infty)$ denote a measurable negative log-likelihood function.
Bayes theorem then yields the posterior probability measure \cite*{stuart2010}
\begin{align}\label{PostMea}
  \mu^y(\diff x) =\frac{1}{Z}\exp(-\Phi_y(x))\mu_0(\diff x),\quad Z:=\int_{\mathbb{R}^n}\exp(-\Phi_y(x))\mu_0(\diff x).
\end{align}
In this context, the posterior measure plays a crucial role in characterizing both the parameter values and their associated uncertainties.
From \eqref{ForMod}, we define the negative log-likelihood
$$
\Phi_y(u):=\frac{1}{2}|\Sigma_{\eta}^{-\frac{1}{2}}(\mathcal{G}(u)-y)|^2.
$$
Combining the prior probability measure $\mu_0$ with density $\pi_0$ and \eqref{PostMea} gives the posterior density up to a normalizing constant
\begin{align}\label{PosDen}
\pi^y(u) = \exp\bigl(-\Phi_y(u)\bigr)\pi_0(u).
\end{align}

\subsection{Evaluation of effective sample size }\label{AppendixESSHD}
Here we describe the calculation of effective sample size (ESS) used in our results.
Let the $\tau_i$ be the integrated autocorrelation time of dimension $i$.
This value is given by
$$
\tau_i=1+2\sum_{j=1}^{N_s}corr(\theta_1,\theta_{1+j})
$$
for dimension $i$ of samples $\{\theta_j\}_{j=1}^{N_s}$, where $corr(\cdot,\cdot)$ is the correlation coefficient.
Then we define the maximum integrated autocorrelation time over all dimension
$$
\tau_{max}=\max_{i\in\{1,2,\dots,n\}}\tau_i.
$$
The ESS is then given by
$$
ESS=\frac{N_s}{\tau_{max}}.
$$

\section{Details about numerical experiments }\label{Appendix:DetailsNumExp}
\subsection{Euler Bernoulli beam problem}\label{Appendix:DetailsNumExpEBb}
In section \ref{subsubsec:EBbProblem}, we consider the Euler Bernoulli beam problem described by PDE \eqref{EBbProblemPDE} for $\Omega=[0,L]$, where $u:\Omega\to\mathbb{R}$ is the vertical deflection of the beam and $f:\Omega\to\mathbb{R}$ is the load.
The beam is in cantilever configuration where the left boundary is fixed and the right boundary is free, i.e., the boundary conditions are
\begin{align}\label{EBbProblemPDECondition}
u(0)=0,\quad \frac{\partial u}{\partial x}\bigg|_{x=0}=0,\quad \frac{\partial^2u}{\partial x}\bigg|_{x=L}=0,\quad \frac{\partial^3u}{\partial x}\bigg|_{x=L}=0.
\end{align}
Consider the function
$$
Z(x,\alpha)=\biggl(1+\exp(-\frac{x-\alpha}{0.005})\biggr)^{-1},
$$
with
$$
\lim_{x\to-\infty}Z(x,\alpha)=0,\quad \lim_{x\to\infty}Z(x,\alpha)=1,
$$
such that there is a smooth transition from $0$ to $1$ at $\alpha$.
Let $L=1$, $\alpha_1=0, \alpha_2=0.1,\alpha_3=1$ and consider the parameter $\theta=(\theta_1,\theta_2)\in\mathbb{R}^2_+$.
Define the function $\widetilde{E}_2:\Omega\times\mathbb{R}\to\mathbb{R}$ as
$$
\widetilde{E}_2(x,\theta)=(1-Z(x,\alpha_2))\theta_1+Z(x,\alpha_2)\theta_2.
$$
Given a parameter $\theta$, the function $\widetilde{E}_2$ is a smooth approximation of the piecewise constant function $\sum_{i=1}^{2}\theta_i\mathbb{I}_{(\alpha_i,\alpha_{i+1}]}$, where $\mathbb{I}_{(\alpha_i,\alpha_{i+1}]}$ is the indicator function of the interval $(\alpha_i,\alpha_{i+1}]\subset \mathbb{R}$.

In the numerical simulation, we use the finite difference to solve the equation \eqref{EBbProblemPDE} and \eqref{EBbProblemPDECondition} on a mesh of $601$ equidistant grid points in $[0,1]$.
We take the load $f(x)=1$.
The Gaussian prior $\mu_0$ is set by $(1,1)$ mean and $[0.3\;0;0\;0.5]^2$ covariance matrix.
The measurement data are obtained by $y=\mathcal{G}(\theta)+\eta$ with $\theta=(1,1)$, where $\mathcal{G}$ is forward model \eqref{EBbProblemPDE} and \eqref{EBbProblemPDECondition} and $\eta$ is  the Gaussian noise with the standard deviation taken by $20\%$ (i.e., noise level) of the last component of the forward model output, i.e., the displacement  $u$ at the end of the beam.
The measurements $y_1,\dots,y_{41}$ are equidistant points of $601$ grid points.
Note that the left end-point at $x = 0$ is not observed since it is fixed by the boundary
condition $u(0)=0$.

\subsection{Locating acoustic sources}\label{Appendix:DetailsNumExpLAS}
For the locating acoustic sources problem, we consider the Helmholtz equation \eqref{Hel}, where $k>0$ is the wave number and the scattered field $u^s$ satisfies the Sommerfeld Radiation Condition
\begin{align}\label{HelSRC}
\lim_{r\rightarrow\infty}\sqrt{r}(\frac{\partial u}{\partial r}-iku)=0,
\end{align}
where $r=|x|,\; i=\sqrt{-1}$, and the sources is given by \eqref{Exa1:source}.
Furthermore, the points $z_i$ are assumed to be mutually distinct.
The solution to the \eqref{Hel} is formally formulated as
$$
u(x)=-\sum_{i=1}^N\lambda_i\Phi(x,z_i),
$$
where
\begin{align}\label{FunSolHel}
\Phi(x,y)=\frac{i}{4}H^{(1)}_0(k|x-y|),
\end{align}
is the fundamental solution to the Helmholtz equation and $H^{(1)}_0$ is Hankel function of the first kind of order $0$.
By the asymptotic behavior of $\Phi$, we obtain the far field pattern
\begin{align}\label{FarFie}
u_{\infty}(\hat{x})=-\frac{e^{i\frac{\pi}{4}}}{\sqrt{8\pi k}}\sum_{i=1}^{N}\lambda_ie^{-ik\hat{x}\cdot z_i},\quad\hat{x}\in S^1.
\end{align}

In this example, we take $\lambda=1$, $k=1$, and the number of direction of is taken by $200$.
The exact location of acoustic source in \eqref{Hel} is set $z_1=(1,3),z_2=(2,2.5),z_3=(1.5,3.5)$ and let $\theta=(1,1.5,2,2.5,3,3.5)$.
The Gaussian prior $\mu_0$ is taken by $[2,2,2,2,2,2]$ mean and  identity covariance matrix.
The measurement data are obtained by $y=\mathcal{G}(\theta)+\eta$ where $\mathcal{G}$ is the far field pattern \eqref{FarFie} and $\eta$ is the Gaussian noise with the standard deviation taken by $1\%,5\%,10\%$ (i.e., noise level) of the maximum norm of the model output.

The target measure for GOAS is the corresponding Bayesian posterior measure. To construct a discrete approximation of the posterior, two density-based discretizations are considered.
For the first design, a preliminary MCMC chain of length $1\times10^4$ is partitioned into $K$ clusters using the $K$-means algorithm \cite{lloyd1982}, and the cluster centers are taken as representative target points. This produces a space-filling, low-discrepancy-like point set concentrated in the dominant posterior region. 

For the second design, we use an importance-sampling strategy. A pilot MCMC ($10^4$ iterations) run is first used only to estimate the posterior mean and covariance, denoted by \(\widehat m\) and \(\widehat C\), respectively. These estimates are used to define a Gaussian proposal
$
q(z)=\mathcal N(z;\widehat m,\widehat C).
$
We then generate \(K\) independent candidate points $z_1,\dots,z_K$ from $q$ and assign the unnormalized importance weights
\[
\widetilde w_i
=
\frac{\pi(z_i)}{q(z_i)},
\]
where \(\widetilde\pi\) denotes the unnormalized posterior density. After normalization,
\[
w_i
=
\frac{\widetilde w_i}
{\sum_{j=1}^K\widetilde w_j},
\]
the posterior is approximated by the weighted discrete measure
\[
\mu_t^K
=
\sum_{i=1}^K w_i\delta_{z_i}.
\]

\subsection{Inverse scattering from an open arc}\label{Appendix:DetailsNumExpIPopenarc}
In section \ref{subsec:exa2}, we consider the inverse scattering from an open arc described by PDE \eqref{IPopenarcPDE}.
Assume that $\Gamma\subset\mathbb{R}^2$ is an arc of class $C^3$, that is
$$
\Gamma:=\{z(s):s\in[-1,1]\},
$$
where $z:[-1,1]\to\mathbb{R}^2$  is an injective and three times continuously differentiable function.
This direct scattering problem is a special case of the exterior Dirichlet problem for the Helmholtz equation.
We can search the scatted field in \eqref{IPopenarcPDE} in the  form of a single-layer potential
\begin{align}\label{SlP}
u^s(x)=\int_{\Gamma}\Phi(x,y)\varphi(y)\diff s(y),\quad x\in\mathbb{R}^2\setminus\Gamma,
\end{align}
where $\Phi$ is the fundamental solution to the Helmholtz equation in two dimensions given by \eqref{FunSolHel}.
 From the potential theoretic jump relations it follows that the potential $u^s$ given by \eqref{SlP} solves the exterior Dirichlet problem provided the density $\varphi$ is a solution to the integral equation \cite{colton1998,colton2013}
\begin{align}\label{exa2IntEqu}
\int_{\Gamma}\Phi(x,y)\varphi(y)\diff s(y)=-u^i(x),\quad x\in\Gamma,
\end{align}
where the incident wave, in the example, is set the plane wave in the direction $\vartheta$, i.e., $u^i(x)=e^{ikx\cdot \vartheta}$.
From the asymptotics for the Hankel function $H^{(1)}_0$, we see that the far-field pattern of the single-layer potential \eqref{SlP} is given by
$$
u_{\infty}(\hat{x})=\frac{e^{i\pi/4}}{\sqrt{8\pi k}}\int_{\Gamma}e^{-ik\hat{x}\cdot y}\varphi(y)\diff s(y),\quad \hat{x}\in S^1.
$$

We use the spectral method to solve the boundary integral equation \eqref{exa2IntEqu}.
The true sound-soft arc $\Gamma$ is described by the parametric representation
\begin{align}\label{arcpara}
z(s)=(\cos s+0.65\cos 2s-0.65,1.5\sin s), \quad \frac{\pi}{4}\leq s\leq \frac{7\pi}{4}.
\end{align}
For simplicity, we denote $z$ as
$$
z(s)=(z_1(s),z_2(s)), \quad \frac{\pi}{4}\leq s\leq \frac{7\pi}{4},
$$
where $z_i(s)=\sum_{j=1}^{N}(c_{i0}+c_{ij}\sin(js)+c'_{ij}\cos(js))$.
Let $N=2$ and $\theta=(c_{10},c_{11},c'_{11},c_{12},c'_{12}, \\ c_{20},c_{21},c'_{21},c_{22},c'_{22})$.
For this inverse problem, we reconstruct the parameter $\theta$.
The true arc $\Gamma$, i.e., \eqref{arcpara} is $(-0.65,0,1,0,0.65,0,1.5,0,0,0)$.
We take $k=1,\vartheta=(1,0)$.
The Gaussian prior is set to $(-0.5,0, 2, 0, 0, 0, 2, 0, 0, 0)$ mean and $0.1^2$ variance.
The measurements are obtained by the far-field pattern at $200$ uniform grid points on the unit circle polluted with zero-mean Gaussian noise with variance $10\%$ of the maximum norm of the far-field pattern output.
Similarly to locating acoustic sources example, we use a length $2\times10^4$ chain from MCMC to obtain a low-discrepancy-like point set and use it to discretize the posterior.

\subsection{Application to a nonlinear advection--diffusion--reaction inverse problem}
\label{Appendix:DetailsNumExpADR}
We apply GOAS to a nonlinear PDE-constrained Bayesian inverse problem involving the simultaneous reconstruction of two unknown parameter fields. Consider the time-dependent advection--diffusion--reaction (ADR) equation
\begin{align}
\label{ADR}
\frac{\partial u}{\partial t}
+\nu\cdot\nabla u
-\nabla\cdot(m_d\nabla u)
+cu^3
&=f,\quad \text{in } \Omega\times(0,T),\\
\frac{\partial u}{\partial \boldsymbol n}
&=0,\quad \text{on } \partial\Omega\times(0,T),\\
u(\cdot,0)
&=m_0,\quad \text{in } \Omega,
\label{ADR3}
\end{align}
where $\boldsymbol n$ denotes the outward unit normal to $\partial\Omega$, $\nu$ is the prescribed advection velocity, $c$ is the reaction coefficient, and $f$ is the source term.
The inverse problem is to determine both the diffusion coefficient filed and initial condition field $(m_d,m_0)$ in ADR model using the measurements data $u(x,T),x\in\Omega$ \cite{ghattas2021learning}.
Owing to the cubic reaction term, the forward map is nonlinear, resulting in a challenging nonlinear Bayesian inverse problem.

In this numerical example, we use the finite element method with Newton iteration to solve the equation \eqref{ADR}-\eqref{ADR3}.
We take $T=0.005,\; \Omega=[0,1]\times [0,1],\; c=1,\; \nu=(cos(t),sin(t))$ and
$$
f(x,t)=\exp\biggl(\frac{\|x-0.5\|_2^2}{0.9^2}\biggr).
$$
The number of time steps of discretization is $21$ and the $\Omega$ is discretised into a triangular mesh with $328$ elements and $185$ vertices.
Let $x=(x_1,x_2)$ and the basis of space $L^2(\Omega)$ is truncated as
$
\{x_1,x_2,x_1^2,x_1x_2,x_2^2,\dots,x_1^{\mathbb{K}},x_1^{\mathbb{K}-1}x_2,\dots,x_2^{\mathbb{K}}\},
$
where $\mathbb{K}$ is a positive integer.
For computational simplicity, in this example we will reconstruct the coefficients of $(m_d,m_0)$ in this polynomial basis.
The exact diffusion coefficient filed and initial condition field in \eqref{ADR} are set
$$
m_0(x)=0.2x_1+0.3x_2+0.4x_1x_2+0.5x_2^2+0.1x_1^3,\quad m_d(x)=0.01x_1,
$$
respectively.
Gaussian prior $\mu_0$ is specified with a mean of $[0,0.15,0.15,0.15,0.15,0.15,0.15]$ and a covariance matrix given by $\text{diag}([0.0001,0.01,0.01,0.01,0.01,0.01,0.01])$.
The measurement data are obtained using the equation $y=\mathcal{G}(m_d,m_0)+\eta$ where $\mathcal{G}$ is the forward model \eqref{ADR}-\eqref{ADR3} and $\eta$ represents the Gaussian noise with a standard deviation taken by $1\%$ (i.e., noise level) of the maximum norm of the model output $u(x,T)$.
Notice that, in order to avoid `inverse crimes', we generate measurement data by solving the forward problem on a finer grid.
An addition, we also use a low-discrepancy-like desugn to discretize the posterior.

\begin{figure}[htbp]
  \centering
  \subfloat[True, $m_d$]
  {
   \includegraphics[width=0.25\textwidth,height=0.2\textwidth]{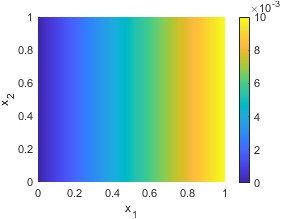}
  }
  \subfloat[GOAS, $m_d$]
  {
    \includegraphics[width=0.25\textwidth,height=0.2\textwidth]{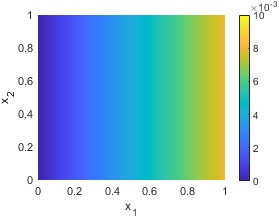}
  }
  \subfloat[Std, $m_d$]
  {
   \includegraphics[width=0.25\textwidth,height=0.2\textwidth]{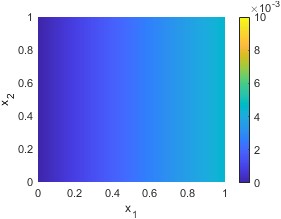}
  }\\
  \subfloat[True, $m_0$]
  {
   \includegraphics[width=0.25\textwidth,height=0.2\textwidth]{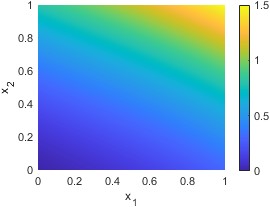}
  }
  \subfloat[GOAS, $m_0$]
  {
   \includegraphics[width=0.25\textwidth,height=0.2\textwidth]{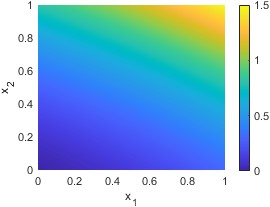}
  }
  \subfloat[Std, $m_0$]
  {
    \includegraphics[width=0.25\textwidth,height=0.2\textwidth]{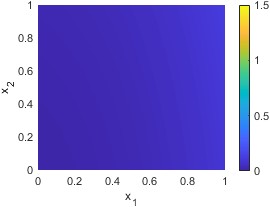}
  }
  \caption{
True fields, posterior sample means, and posterior sample standard deviations (Std) obtained by GOAS with $K=350$ for the simultaneous reconstruction of the diffusion coefficient field $m_d$ (top row) and the initial condition field $m_0$ (bottom row) in the nonlinear advection--diffusion--reaction inverse problem.}
\label{fig:adr}
\end{figure}

Figure~\ref{fig:adr} shows the posterior sample means and standard deviations obtained by GOAS with $K=350$.
The posterior mean fields closely reproduce the main spatial features of both the diffusion coefficient $m_d$ and the initial condition $m_0$. The corresponding posterior standard deviations are relatively small, indicating that the observations provide substantial information about both unknown fields.

This example illustrates that the GOAS framework can also be applied to nonlinear PDE-constrained Bayesian inverse problems involving multiple unknown parameter fields.


\section{Proofs}

\begin{lem}\label{lem:Treg}
Let $\mathcal{M}_1$ and $\mathcal{M}_2$ be two compact manifold in $\mathbb{R}^{n+1}$ endowed with two geodesic distance $d_{\mathcal{M}_1}$ and $d_{\mathcal{M}_2}$, respectively.
Let the $\mu_1$ and $\mu_2$ be two measures on the manifold $\mathcal{M}_1$.
If the map $T:\mathcal{M}_1\to \mathcal{M}_2$ is Lipschitz continuous, then
$$
W_p(T_{\sharp}\mu_1,T_{\sharp}\mu_2)\leq LW_p(\mu_1,\mu_2),
$$
where $L$ is the Lipschitz constant of $T$.
\end{lem}
\begin{proof}
Let $\gamma\in\Gamma(u_1,u_2)$, and define a measure
$$
\gamma_T:=(T\times T)_{\sharp}\gamma,
$$
which is a coupling of $T_{\sharp}u_1$ and $T_{\sharp}u_2$, i.e., $\gamma_T(B\times C)=\gamma(T^{-1}(B)\times T^{-1}(C))$ for any Borel set $B,C\subset\mathcal{M}_2$.
Then
\begin{align}\label{error:eq2}
W_p(T_{\sharp}u_1,T_{\sharp}u_2)&=
\biggl(\inf_{\gamma_T\in\Gamma(T_{\sharp}u_1,T_{\sharp}u_2)}\int_{\mathcal{M}_2\times \mathcal{M}_2}d^p_{\mathcal{M}_2}(x,y)\diff\gamma_T(x,y)\biggr)^{1/p}\notag\\
&\leq\biggl(\int_{\mathcal{M}_2\times \mathcal{M}_2}d^p_{\mathcal{M}_2}(x,y)\diff\gamma_T(x,y)  \biggr)^{1/p}\notag \\
&=\biggl(\int_{\mathcal{M}_1\times \mathcal{M}_1}d^p_{\mathcal{M}_2}(T(x),T(y))\diff\gamma(x,y)  \biggr)^{1/p}.
\end{align}
Since the map $T:\mathcal{M}_1\to \mathcal{M}_2$ is Lipschitz continuous, we have
$$
d^p_{\mathcal{M}_2}(T(x),T(y))\leq L^pd^p_{\mathcal{M}_1}(x,y),
$$
where $L$ is the Lipschitz constant of $T$.
Plugging this bound into \eqref{error:eq2}, we deduce that
$$
W_p(T_{\sharp}u_1,T_{\sharp}u_2)\leq \biggl(\int_{\mathcal{M}_1\times \mathcal{M}_1}L^pd^p_{\mathcal{M}_1}(x,y)\diff\gamma(x,y)  \biggr)^{1/p},
$$
which implies
$$
W_p(T_{\sharp}\mu_1,T_{\sharp}\mu_2)\leq LW_p(\mu_1,\mu_2).
$$
\end{proof}

\end{appendices}




\bibliographystyle{siamplain}
 \bibliography{bibliography.bib}

\begin{thebibliography}{10}

\bibitem{andrieu2006}
{\sc C.~Andrieu and {\'E}.~Moulines}, {\em On the ergodicity properties of some
  adaptive {MCMC} algorithms}, The Annals of Applied Probability, 16 (2006),
  pp.~1462--1505, \url{https://doi.org/10.1214/105051606000000286}.

\bibitem{baptista2023}
{\sc R.~Baptista, Y.~Marzouk, and O.~Zahm}, {\em On the representation and
  learning of monotone triangular transport maps}, Foundations of Computational
  Mathematics, 24 (2024), pp.~2063--2108,
  \url{https://doi.org/10.1007/s10208-023-09630-x}.

\bibitem{benamou2020}
{\sc J.-D. Benamou, W.~L. Ijzerman, and G.~Rukhaia}, {\em An entropic optimal
  transport numerical approach to the reflector problem}, Methods and
  Applications of Analysis, 27 (2020), pp.~311--340,
  \url{https://doi.org/10.4310/MAA.2020.v27.n4.a1}.

\bibitem{steve2011}
{\sc S.~Brooks, A.~Gelman, G.~L. Jones, and X.-L. Meng}, eds., {\em Handbook of
  Markov Chain Monte Carlo}, Chapman and Hall/CRC, Boca Raton, FL, 2011,
  \url{https://doi.org/10.1201/b10905}.

\bibitem{caffarelli1999}
{\sc L.~A. Caffarelli, S.~A. Kochengin, and V.~I. Oliker}, {\em On the
  numerical solution of the problem of reflector design with given far-field
  scattering data}, in Monge Amp\`ere Equation: Applications to Geometry and
  Optimization, L.~A. Caffarelli and M.~Milman, eds., vol.~226 of Contemporary
  Mathematics, American Mathematical Society, Providence, RI, 1999, pp.~13--32,
  \url{https://doi.org/10.1090/conm/226/03233}.

\bibitem{caffarelli2008}
{\sc L.~A. Caffarelli and V.~I. Oliker}, {\em Weak solutions of one inverse
  problem in geometric optics}, Journal of Mathematical Sciences, 154 (2008),
  pp.~39--49, \url{https://doi.org/10.1007/s10958-008-9152-x}.

\bibitem{colton1998}
{\sc D.~Colton and R.~Kress}, {\em Inverse Acoustic and Electromagnetic
  Scattering Theory}, vol.~93 of Applied Mathematical Sciences,
  Springer-Verlag, Berlin, 2~ed., 1998,
  \url{https://doi.org/10.1007/978-3-662-03537-5}.

\bibitem{colton2013}
{\sc D.~Colton and R.~Kress}, {\em Integral Equation Methods in Scattering
  Theory}, vol.~72 of Classics in Applied Mathematics, Society for Industrial
  and Applied Mathematics, Philadelphia, PA, 2013,
  \url{https://doi.org/10.1137/1.9781611973167}.

\bibitem{cotter2013}
{\sc S.~L. Cotter, G.~O. Roberts, A.~M. Stuart, and D.~White}, {\em Mcmc
  methods for functions: Modifying old algorithms to make them faster},
  Statistical Science, 28 (2013), pp.~424--446,
  \url{https://doi.org/10.1214/13-STS421}.

\bibitem{cui2016}
{\sc T.~Cui, K.~J.~H. Law, and Y.~M. Marzouk}, {\em Dimension-independent
  likelihood-informed {MCMC}}, Journal of Computational Physics, 304 (2016),
  pp.~109--137, \url{https://doi.org/10.1016/j.jcp.2015.10.008}.

\bibitem{duane1987}
{\sc S.~Duane, A.~D. Kennedy, B.~J. Pendleton, and D.~Roweth}, {\em Hybrid
  monte carlo}, Physics Letters B, 195 (1987), pp.~216--222,
  \url{https://doi.org/10.1016/0370-2693(87)91197-X}.

\bibitem{durkan2019}
{\sc C.~Durkan, A.~Bekasov, I.~Murray, and G.~Papamakarios}, {\em Neural spline
  flows}, in Advances in Neural Information Processing Systems, vol.~32, 2019,
  pp.~7509--7520.

\bibitem{el2011}
{\sc A.~El~Badia and T.~Nara}, {\em An inverse source problem for {Helmholtz}'s
  equation from the {Cauchy} data with a single wave number}, Inverse Problems,
  27 (2011), p.~105001, \url{https://doi.org/10.1088/0266-5611/27/10/105001}.

\bibitem{el2012}
{\sc T.~A. El~Moselhy and Y.~M. Marzouk}, {\em Bayesian inference with optimal
  maps}, Journal of Computational Physics, 231 (2012), pp.~7815--7850,
  \url{https://doi.org/10.1016/j.jcp.2012.07.022}.

\bibitem{eller2009}
{\sc M.~Eller and N.~P. Valdivia}, {\em Acoustic source identification using
  multiple frequency information}, Inverse Problems, 25 (2009), p.~115005,
  \url{https://doi.org/10.1088/0266-5611/25/11/115005}.

\bibitem{fournier2010freeform}
{\sc F.~Fournier}, {\em Freeform Reflector Design with Extended Sources}, PhD
  thesis, University of Central Florida, Orlando, FL, 2010,
  \url{https://stars.library.ucf.edu/etd/4283/}.

\bibitem{fournier2015}
{\sc N.~Fournier and A.~Guillin}, {\em On the rate of convergence in
  wasserstein distance of the empirical measure}, Probability Theory and
  Related Fields, 162 (2015), pp.~707--738,
  \url{https://doi.org/10.1007/s00440-014-0583-7}.

\bibitem{ghattas2021learning}
{\sc O.~Ghattas and K.~Willcox}, {\em Learning physics-based models from data:
  Perspectives from inverse problems and model reduction}, Acta Numerica, 30
  (2021), pp.~445--554, \url{https://doi.org/10.1017/S0962492921000064}.

\bibitem{girolami2011}
{\sc M.~Girolami and B.~Calderhead}, {\em Riemann manifold {Langevin} and
  {Hamiltonian} monte carlo methods}, Journal of the Royal Statistical Society:
  Series B (Statistical Methodology), 73 (2011), pp.~123--214,
  \url{https://doi.org/10.1111/j.1467-9868.2010.00765.x}.

\bibitem{graf2012}
{\sc T.~Graf and V.~I. Oliker}, {\em An optimal mass transport approach to the
  near-field reflector problem in optical design}, Inverse Problems, 28 (2012),
  p.~025001, \url{https://doi.org/10.1088/0266-5611/28/2/025001}.

\bibitem{grathwohl2019}
{\sc W.~Grathwohl, R.~T.~Q. Chen, J.~Bettencourt, I.~Sutskever, and
  D.~Duvenaud}, {\em {FFJORD}: Free-form continuous dynamics for scalable
  reversible generative models}, in International Conference on Learning
  Representations, 2019, \url{https://openreview.net/forum?id=rJxgknCcK7}.

\bibitem{guan1998}
{\sc P.~Guan and X.-J. Wang}, {\em On a {Monge--Amp\`ere} equation arising in
  geometric optics}, Journal of Differential Geometry, 48 (1998), pp.~205--223,
  \url{https://doi.org/10.4310/jdg/1214460795}.

\bibitem{haario2001}
{\sc H.~Haario, E.~Saksman, and J.~Tamminen}, {\em An adaptive metropolis
  algorithm}, Bernoulli, 7 (2001), pp.~223--242,
  \url{https://doi.org/10.2307/3318737}.

\bibitem{jaini2019}
{\sc P.~Jaini, K.~A. Selby, and Y.~Yu}, {\em Sum-of-squares polynomial flow},
  in Proceedings of the 36th International Conference on Machine Learning,
  vol.~97 of Proceedings of Machine Learning Research, PMLR, 2019,
  pp.~3009--3018, \url{https://proceedings.mlr.press/v97/jaini19a.html}.

\bibitem{karakhanyan2010}
{\sc A.~Karakhanyan and X.-J. Wang}, {\em On the reflector shape design},
  Journal of Differential Geometry, 84 (2010), pp.~561--610,
  \url{https://doi.org/10.4310/jdg/1279114301}.

\bibitem{kirsch2000}
{\sc A.~Kirsch and S.~Ritter}, {\em A linear sampling method for inverse
  scattering from an open arc}, Inverse Problems, 16 (2000), pp.~89--105,
  \url{https://doi.org/10.1088/0266-5611/16/1/308}.

\bibitem{kobyzev2020}
{\sc I.~Kobyzev, S.~J.~D. Prince, and M.~A. Brubaker}, {\em Normalizing flows:
  An introduction and review of current methods}, IEEE Transactions on Pattern
  Analysis and Machine Intelligence, 43 (2021), pp.~3964--3979,
  \url{https://doi.org/10.1109/TPAMI.2020.2992934}.

\bibitem{kochengin2003}
{\sc S.~A. Kochengin and V.~I. Oliker}, {\em Computational algorithms for
  constructing reflectors}, Computing and Visualization in Science, 6 (2003),
  pp.~15--21, \url{https://doi.org/10.1007/s00791-003-0103-2}.

\bibitem{kress1995}
{\sc R.~Kress}, {\em Inverse scattering from an open arc}, Mathematical Methods
  in the Applied Sciences, 18 (1995), pp.~267--293,
  \url{https://doi.org/10.1002/mma.1670180403}.

\bibitem{lipman2023}
{\sc Y.~Lipman, R.~T.~Q. Chen, H.~Ben-Hamu, M.~Nickel, and M.~Le}, {\em Flow
  matching for generative modeling}, in International Conference on Learning
  Representations, 2023, \url{https://openreview.net/forum?id=PqvMRDCJT9t}.

\bibitem{liu2021}
{\sc Y.~Liu, Y.~Guo, and J.~Sun}, {\em A deterministic-statistical approach to
  reconstruct moving sources using sparse partial data}, Inverse Problems, 37
  (2021), p.~065005, \url{https://doi.org/10.1088/1361-6420/abf813}.

\bibitem{lloyd1982}
{\sc S.~P. Lloyd}, {\em Least squares quantization in {PCM}}, IEEE Transactions
  on Information Theory, 28 (1982), pp.~129--137,
  \url{https://doi.org/10.1109/TIT.1982.1056489}.

\bibitem{martin2012}
{\sc J.~Martin, L.~C. Wilcox, C.~Burstedde, and O.~Ghattas}, {\em A stochastic
  {Newton} {MCMC} method for large-scale statistical inverse problems with
  application to seismic inversion}, SIAM Journal on Scientific Computing, 34
  (2012), pp.~A1460--A1487, \url{https://doi.org/10.1137/110845598}.

\bibitem{marzouk2016}
{\sc Y.~Marzouk, T.~Moselhy, M.~Parno, and A.~Spantini}, {\em Sampling via
  measure transport: An introduction}, in Handbook of Uncertainty
  Quantification, R.~Ghanem, D.~Higdon, and H.~Owhadi, eds., Springer
  International Publishing, Cham, 2017, pp.~785--825,
  \url{https://doi.org/10.1007/978-3-319-12385-1_23}.

\bibitem{neal2003}
{\sc R.~M. Neal}, {\em Slice sampling}, The Annals of Statistics, 31 (2003),
  pp.~705--767, \url{https://doi.org/10.1214/aos/1056562461}.

\bibitem{neal2011}
{\sc R.~M. Neal}, {\em {MCMC} using {Hamiltonian} dynamics}, in Handbook of
  Markov Chain Monte Carlo, S.~Brooks, A.~Gelman, G.~L. Jones, and X.-L. Meng,
  eds., Chapman and Hall/CRC, Boca Raton, FL, 2011, pp.~113--162,
  \url{https://doi.org/10.1201/b10905-6}.

\bibitem{oliker1993}
{\sc V.~I. Oliker and E.~Newman}, {\em The energy conservation equation in the
  reflector mapping problem}, Applied Mathematics Letters, 6 (1993),
  pp.~91--95, \url{https://doi.org/10.1016/0893-9659(93)90156-H}.

\bibitem{papamakarios2021}
{\sc G.~Papamakarios, E.~Nalisnick, D.~J. Rezende, S.~Mohamed, and
  B.~Lakshminarayanan}, {\em Normalizing flows for probabilistic modeling and
  inference}, Journal of Machine Learning Research, 22 (2021), pp.~1--64,
  \url{https://jmlr.org/papers/v22/19-1028.html}.

\bibitem{parno2022}
{\sc M.~Parno, P.-B. Rubio, D.~Sharp, M.~Brennan, R.~Baptista, H.~Bonart, and
  Y.~Marzouk}, {\em {MParT}: Monotone parameterization toolkit}, Journal of
  Open Source Software, 7 (2022), p.~4843,
  \url{https://doi.org/10.21105/joss.04843}.

\bibitem{parno2018}
{\sc M.~D. Parno and Y.~M. Marzouk}, {\em Transport map accelerated markov
  chain monte carlo}, SIAM/ASA Journal on Uncertainty Quantification, 6 (2018),
  pp.~645--682, \url{https://doi.org/10.1137/17M1134640}.

\bibitem{peherstorfer2019}
{\sc B.~Peherstorfer and Y.~Marzouk}, {\em A transport-based multifidelity
  preconditioner for markov chain monte carlo}, Advances in Computational
  Mathematics, 45 (2019), pp.~2321--2348,
  \url{https://doi.org/10.1007/s10444-019-09711-y}.

\bibitem{peyre2019computational}
{\sc G.~Peyr{\'e} and M.~Cuturi}, {\em Computational optimal transport: With
  applications to data science}, Foundations and Trends in Machine Learning, 11
  (2019), pp.~355--607, \url{https://doi.org/10.1561/2200000073}.

\bibitem{rezende2015}
{\sc D.~J. Rezende and S.~Mohamed}, {\em Variational inference with normalizing
  flows}, in Proceedings of the 32nd International Conference on Machine
  Learning, vol.~37 of Proceedings of Machine Learning Research, PMLR, 2015,
  pp.~1530--1538, \url{https://proceedings.mlr.press/v37/rezende15.html}.

\bibitem{robert2004}
{\sc C.~P. Robert and G.~Casella}, {\em Monte Carlo Statistical Methods},
  Springer Texts in Statistics, Springer, New York, 2~ed., 2004,
  \url{https://doi.org/10.1007/978-1-4757-4145-2}.

\bibitem{roberts1998}
{\sc G.~O. Roberts and J.~S. Rosenthal}, {\em Optimal scaling of discrete
  approximations to {Langevin} diffusions}, Journal of the Royal Statistical
  Society: Series B (Statistical Methodology), 60 (1998), pp.~255--268,
  \url{https://doi.org/10.1111/1467-9868.00123}.

\bibitem{shirley2008realistic}
{\sc P.~Shirley and R.~K. Morley}, {\em Realistic Ray Tracing}, A K Peters,
  Natick, MA, 2~ed., 2008.

\bibitem{stuart2010}
{\sc A.~M. Stuart}, {\em Inverse problems: A {Bayesian} perspective}, Acta
  Numerica, 19 (2010), pp.~451--559,
  \url{https://doi.org/10.1017/S0962492910000061}.

\bibitem{sullivan2010}
{\sc A.~B. Sullivan, D.~M. Snyder, and S.~A. Rounds}, {\em Controls on
  biochemical oxygen demand in the upper {Klamath River}, {Oregon}}, Chemical
  Geology, 269 (2010), pp.~12--21,
  \url{https://doi.org/10.1016/j.chemgeo.2009.08.007}.

\bibitem{sun2024}
{\sc Z.~Sun and G.-H. Zheng}, {\em Geometric optics approximation sampling:
  Near-field case}, arXiv.2403.01655,  (2024),
  \url{https://doi.org/10.48550/arXiv.2403.01655}.

\bibitem{tierney1994}
{\sc L.~Tierney}, {\em Markov chains for exploring posterior distributions},
  The Annals of Statistics, 22 (1994), pp.~1701--1728,
  \url{https://doi.org/10.1214/aos/1176325750}.

\bibitem{villani2021}
{\sc C.~Villani}, {\em Topics in Optimal Transportation}, vol.~58 of Graduate
  Studies in Mathematics, American Mathematical Society, Providence, RI, 2003.

\bibitem{villani2009}
{\sc C.~Villani}, {\em Optimal Transport: Old and New}, vol.~338 of Grundlehren
  der Mathematischen Wissenschaften, Springer, Berlin, Heidelberg, 2009,
  \url{https://doi.org/10.1007/978-3-540-71050-9}.

\bibitem{wang1996}
{\sc X.-J. Wang}, {\em On the design of a reflector antenna}, Inverse Problems,
  12 (1996), pp.~351--375, \url{https://doi.org/10.1088/0266-5611/12/3/013}.

\bibitem{wang2004}
{\sc X.-J. Wang}, {\em On the design of a reflector antenna {II}}, Calculus of
  Variations and Partial Differential Equations, 20 (2004), pp.~329--341,
  \url{https://doi.org/10.1007/s00526-003-0239-4}.

\bibitem{wenliang2019}
{\sc L.~Wenliang, D.~J. Sutherland, H.~Strathmann, and A.~Gretton}, {\em
  Learning deep kernels for exponential family densities}, in Proceedings of
  the 36th International Conference on Machine Learning, vol.~97 of Proceedings
  of Machine Learning Research, PMLR, 2019, pp.~6737--6746,
  \url{https://proceedings.mlr.press/v97/wenliang19a.html}.

\end{thebibliography}

\end{document}